\documentclass[a4paper,11pt,twoside]{article}
\usepackage{pst-fill,pst-grad}
\usepackage{textcomp}
\usepackage[english]{babel}
\usepackage{graphicx}
\usepackage{amsmath}
\usepackage{float}
\usepackage[matrix,arrow,curve]{xy}
\usepackage{pstricks}
\usepackage{mathtools}
\usepackage{amsmath,amsfonts,verbatim,afterpage,theorem,euscript,mathrsfs,amssymb}
\usepackage{amsfonts}
\usepackage{amssymb}
\usepackage{array}
\usepackage{dsfont}
\usepackage{hyperref}
\def \vu{\vec{u}}
\def \vpsi{\vec{\psi}}
\def \vf{\vec{f}}
\def \vg{\vec{g}}
\def \div{\mathrm{div}}
\def \vn{\vec{\nabla}}

\newtheorem{Theoreme}{Theorem}

\newtheorem{Proposition}{Proposition}[section]
\newtheorem{Lemme}{Lemma}[section]

\newtheorem{Remarque}{Remark}[section]

\numberwithin{equation}{section}
\title{\bf Some general external forces and critical mild solutions for the fractional Navier-Stokes equations}
\usepackage{authblk}
\author{Diego Chamorro\footnote{\emph{diego.chamorro@univ-evry.fr} }}
\affil{\footnotesize LaMME, Univ. Evry, CNRS, Universit\'e Paris-Saclay, 91025, Evry, France.}
\author{Maxence Mansais\footnote{\emph{maxence.mansais@ens-paris-saclay.fr} }}
\affil{\footnotesize LaMME, Univ. Evry, CNRS, Universit\'e Paris-Saclay, 91025, Evry, France.}
\begin{document}
\maketitle
\begin{abstract}
In this article we study mild solutions for the forced, incompressible fractional Navier-Stokes equations. These solutions are classically obtained via a fixed-point argument which relies on suitable estimates for the initial data, the nonlinearity and the external forces.  Many functional spaces can be considered,  however we are mainly interested here in a critical setting which ensures the existence of global solutions. We give some examples of such critical functional spaces and we discuss their relationship with generic external forces. 
\end{abstract}
\section{The fractional Navier-Stokes equations}

We consider here the fractional incompressible Navier-Stokes equations over $\mathbb{R}^d$ with $d\geq 3$,
\begin{equation}\label{EquationNS_intro}
\begin{cases}
\partial_t \vec{u}  = -(-\Delta)^\frac{\alpha}{2} \vec{u}-\mathrm{div} \left(\vec{u}\otimes\vec{u}\right) -\vn p + \vec{f}, \qquad \mathrm{div}(\vec{u})=0, \qquad (1<\alpha<2),\\[3mm]
\vec{u}(0,.)  = \vec{u_0},\qquad \mathrm{div}(\vec{u_0})=0,
\end{cases}
\end{equation}
where $\vec{u_0}:\mathbb{R}^d\longrightarrow \mathbb{R}^d$ denotes the initial velocity of the fluid, $\vec{u}:[0,+\infty[\times \mathbb{R}^d\longrightarrow \mathbb{R}^d$ is the velocity field which is assumed to be divergence free, $p:[0,+\infty[\times \mathbb{R}^d\longrightarrow \mathbb{R}$ is the internal pressure and $\vec{f}:[0,+\infty[\times \mathbb{R}^d\longrightarrow \mathbb{R}^d$ is a given external force. We recall that for a index $1< \alpha < 2$, the fractional laplacian operator $(-\Delta)^\frac{\alpha}{2}$ can be defined as a Fourier multiplier through the expression $\widehat{(-\Delta)^\frac{\alpha}{2}\phi}(\xi)= |\xi|^\alpha \widehat{\phi}(\xi)$, where $\phi:\mathbb{R}^d\longrightarrow \mathbb{R}$ is a suitable function (say $\phi\in \mathcal{S}(\mathbb{R}^d)$). \\

This fractional system has been studied from different points of view in \cite{LiZhai}, \cite{Miao}, \cite{YuZhai} and \cite{Zhai}, see also the references therein for more details. Let us point out that, just as for the classical Navier-Stokes system, the equation (\ref{EquationNS_intro}) still has many open problems related to existence, regularity and uniqueness issues.\\

We recall now that it is possible to separate the study of the pressure $p$ from the velocity field $\vu$: indeed, by formally applying the divergence to the equation (\ref{EquationNS_intro}) and using the divergence free property of the velocity field ($\div(\vu)=0$) we obtain the following expression
$$(-\Delta) p=\mathrm{div}\left(\mathrm{div}(\vec{u}\otimes \vec{u})-\vec{f}\right),$$
which allows us to deduce information over the pressure as long as we have suitable informations over the velocity field and the external force. To continue, we consider now the Leray projector $\mathbb{P}(\vpsi)=\vpsi + \vn (-\Delta)^{-1}\div(\vpsi)$ and since we have $\mathbb{P}(\vn p)=0$ and $\mathbb{P}(\vu)=\vu$ (since $\div(\vu)=0$), we obtain from (\ref{EquationNS_intro}) the equation
$$\partial_t \vec{u}  = -(-\Delta)^\frac{\alpha}{2} \vec{u}-\mathbb{P}\left(\mathrm{div} \left(\vec{u}\otimes\vec{u}\right)\right)+ \mathbb{P}(\vec{f}),$$
and we consider the integral formulation of the previous equation
\begin{equation}\label{Formulation_Integrale}
\vec{u}=\mathfrak{p}_t\ast \vec{u_0}-\int_0^t \mathfrak{p}_{t-s}\ast\mathbb{P}\left(\mathrm{div}(\vec{u}\otimes \vec{u})\right)ds+\int_0^t \mathfrak{p}_{t-s}\ast\mathbb{P}(\vec{f})ds, 
\end{equation}
where $\mathfrak{p}_t$, with $t>0$, stands for the fractional heat kernel associated to the semi-group $e^{-t(-\Delta)^{\frac{\alpha}{2}}}$ whose action in the Fourier level is given by the formula $\left[e^{-t(-\Delta)^{\frac{\alpha}{2}}} \vpsi\right]^{\widehat{} }(\xi)=e^{-t|\xi|^\alpha}\widehat{\vpsi}(\xi)$: we thus have 
\begin{equation}\label{Definition_FracHeat}
\widehat{\mathfrak{p}_t}(\xi)=e^{-t|\xi|^\alpha}.
\end{equation}
For more details on the properties of this semi-group and its kernel see \cite{kolokoltsov} and \cite{Pazy}.\\

In this article we are mainly interested in the relationships between the initial data $\vu_0$, the general external forces $\vf$ and \emph{mild} solutions $\vu$ of the system (\ref{Formulation_Integrale}) that can be obtained via a fixed point argument in \emph{resolution spaces} that are \emph{critical}. We point out that the main interest for working with \emph{critical spaces} relies on the fact that it is possible to obtain \emph{global} in time solutions.\\ 

Recall that a functional space $(\mathcal{X}, \|\cdot\|_{\mathcal{X}})\subset \mathcal{S}'([0,+\infty[\times \mathbb{R}^d, \mathbb{R}^d)$ is a \emph{resolution space} for the problem (\ref{Formulation_Integrale}) associated to an initial data $\vu_0\in (\mathcal{X}_0, \|\cdot\|_{\mathcal{X}_0})\subset  \mathcal{S}'(\mathbb{R}^d, \mathbb{R}^d)$ and to an external force $\vf \in (\mathcal{Y}, \|\cdot\|_{\mathcal{Y}})\subset \mathcal{S}'([0,+\infty[\times \mathbb{R}^d, \mathbb{R}^d)$ if we have the three following estimates: 
\begin{eqnarray}
&&\|\mathfrak{p}_t\ast \vu_0\|_{\mathcal{X}}\leq C_1\|\vu_0\|_{\mathcal{X}_0},\label{Controle_DonneeInitiale}\\[3mm]
&&\left\|\int_0^t \mathfrak{p}_{t-s}\ast\mathbb{P}\left(\mathrm{div}(\vec{u}\otimes \vec{u})\right)ds\right\|_{\mathcal{X}}\leq C_2\|\vu\|_{\mathcal{X}}\|\vu\|_{\mathcal{X}}\qquad \mbox{and}\label{Controle_Bilineaire}\\[3mm]
&&\left\|\int_0^t \mathfrak{p}_{t-s}\ast\mathbb{P}(\vec{f})ds\right\|_{\mathcal{X}}\leq C_3\|\vf\|_{\mathcal{Y}}.\label{Controle_Force}
\end{eqnarray}
With these estimates at hand (and under a smallness assumption of the quantities $\|\vu_0\|_{\mathcal{X}_0}$ and $\|\vf\|_{\mathcal{Y}}$), it is easy to construct \emph{mild} solutions for the problem (\ref{Formulation_Integrale}) via a Banach-Picard contraction principle, see \cite[Théorème 4.1.1]{Chamorro_Livre}. Many different functional spaces $\mathcal{X}_0$, $\mathcal{X}$ and $\mathcal{Y}$ that satisfy the previous estimates can be considered to solve the integral problem (\ref{Formulation_Integrale}) but we are interested here in a \emph{critical} setting:  indeed, if $\vu(t,x)$ and $\vf(t,x)$ satisfy the equation (\ref{Formulation_Integrale}), then for all $\lambda>0$ the rescaled functions 
\begin{equation}\label{Homogeneite}
\vu_\lambda(t,x)=\lambda^{\alpha-1}\vu(\lambda^\alpha t,\lambda x),
\end{equation}
and
\begin{equation}\label{HomogeneiteForce}
\vf_\lambda(t,x)=\lambda^{2\alpha-1}\vf(\lambda^\alpha t,\lambda x), 
\end{equation}
also satisfy the equation (\ref{Formulation_Integrale}). We will thus say that the functional spaces $(\mathcal{X}, \|\cdot\|_{\mathcal{X}})$ and $(\mathcal{Y}, \|\cdot\|_{\mathcal{Y}})$ are \emph{critical} spaces for the velocity field $\vu$ and the external force $\vf$ if these functional spaces are translation invariant and if moreover, for all $\lambda >0$, we have the homogeneity identities
\begin{equation}\label{Espaces_Critiques}
\|\vu_\lambda\|_{\mathcal{X}}=\|\vu\|_{\mathcal{X}},
\end{equation}
and
\begin{equation}\label{Espaces_CritiquesForce}
\|\vf_\lambda\|_{\mathcal{Y}}=\|\vf\|_{\mathcal{Y}}.
\end{equation}
In the same spirit, we will say that $(\mathcal{X}_0, \|\cdot\|_{\mathcal{X}_0})$ is a \emph{critical} space for the initial data $\vu_0$ if, setting $\vu_{0_\lambda}(x)=\lambda^{\alpha-1}\vu_{0}(\lambda x)$ for all $\lambda>0$, we have the identity 
\begin{equation}\label{Espaces_Critiques_DI}
\|\vu_{0_\lambda}\|_{\mathcal{X}_0}=\|\vu_0\|_{\mathcal{X}_0}.\\[5mm]
\end{equation}

Let us note now that, in the case of the \emph{unforced} classical 3D Navier-Stokes equations (\emph{i.e.} when $\alpha=2$ and $\vf=0$), many different critical spaces have been studied. For example in \cite{FujitaKato} the space $\mathcal{X}=L^\infty_t\dot{H}^{\frac{1}{2}}_x\cap L^4_t\dot{H}^{1}_x$ was studied with an initial data in the space $\mathcal{X}_0=\dot{H}^{\frac{1}{2}}$, in \cite{Kozono} the functional setting $\mathcal{X}=L^\infty_tL^{3,\infty}_x$ was used with $\mathcal{X}_0=L^{3,\infty}$ (where is $L^{3,\infty}$ a Lorentz space) and in \cite{Kato} it was considered $\mathcal{X}=L^\infty_t\mathcal{M}^{2,3}_x\cap \{\underset{t>0}{\sup}\sqrt{t}\|\vu(t, \cdot)\|_{L^\infty} <+\infty\}$ with $\mathcal{X}_0=\mathcal{M}^{2,3}$ (where $\mathcal{M}^{2,3}$ is a Morrey space). Finally, in \cite{Koch}, Koch and Tataru reached what it seems to be -to the best of our knowledge- the largest space where to search for mild solutions, in particular they give a characterization of the corresponding space for the initial data and we have $\mathcal{X}_0=BMO^{-1}$. See the book \cite{PGLR} for a detailed presentation of these (and others) critical spaces used to study the classical Navier-Stokes equation.\\ 

An important remark is the next one: in the previous examples, we have the following relationships for the spaces $\mathcal{X}_0$ used for the initial data $\vu_0$
\begin{equation}\label{Space_Inclusion_X0}
\dot{H}^{\frac{1}{2}}(\mathbb{R}^3)\subset L^{3,\infty}(\mathbb{R}^3)\subset \mathcal{M}^{2,3}(\mathbb{R}^3)\subset BMO^{-1}(\mathbb{R}^3). 
\end{equation}
Note also that we have the space inclusion 
\begin{equation}\label{Space_Inclusion_BMO_Besov}
BMO^{-1}(\mathbb{R}^3)\subset \dot{B}^{-1,\infty}_\infty(\mathbb{R}^3), 
\end{equation}
but as it was observed in \cite{besov_ill_posed}, although the Besov space $\dot{B}^{-1,\infty}_\infty(\mathbb{R}^3)$ is a critical space for the initial data, the \emph{mild} approach is ill-posed in this setting. Remark however that the space $\dot{B}^{-1,\infty}_\infty(\mathbb{R}^3)$ plays an important role in the chain of inclusions (\ref{Space_Inclusion_X0})-(\ref{Space_Inclusion_BMO_Besov}) as, under the translation invariance and homogeneity constraints, it is known since the work \cite{Meyer} that this Besov space is maximal and it thus represents an end point for the previous chain of inclusions.\\

In the case of the \emph{forced} Navier-Stokes equations (\emph{i.e.} with $\vf\neq 0$), in \cite{Cannone} and \cite{Fabes} some functional spaces $\mathcal{Y}$ were considered for the external force. However, as it was recently pointed out in \cite{mixed_norm_time_space}, although some of the properties of the space $\mathcal{X}$ which allow to close the fixed-point argument are ``naturally'' given by the structure of the equation (\emph{i.e.} we must have the estimate (\ref{Controle_Bilineaire})), the choice of the functional space $\mathcal{Y}$ for the external force can be slightly more subtle to study.\\

Coming back to the fractional integral problem (\ref{Formulation_Integrale}) without external force, an adaptation of the Koch and Tataru theorem was made in the article \cite{Zhai} with an initial data $\vu_0$ in the critical space $BMO^{-(\alpha-1)}$. However, and quite surprisingly, due to the properties of the  fractional heat kernel present in the equation (\ref{Formulation_Integrale}), it was proven in \cite{YuZhai} that (and unlike to the case of the classical Navier-Stokes equation) it is possible to close the fixed-point argument for an initial data $\vu_0$ in the maximal critical space $\dot{B}^{-(\alpha-1),\infty}_\infty$ (we have the space inclusion $BMO^{-(\alpha-1)}(\mathbb{R}^d)\subset \dot{B}^{-(\alpha-1),\infty}_\infty(\mathbb{R}^d)$). \\

Since (to the best of our knowledge) the study of the \emph{forced} fractional Navier-Stokes equation was not considered before, we will study here the relationships between external forces $\vf$ and \emph{mild} solutions $\vu$ of the system (\ref{Formulation_Integrale}). Our first result reads as follows:
\begin{Theoreme}\label{Theorem1}
Let $1<\alpha<2$ be fixed. We fix a parameter $p_0$ such that $\frac{d}{\alpha}<p_0 \leq +\infty$ as well as a parameter $\beta>0$ such that $\alpha-\frac{d}{p_0}-1< \beta < \alpha-\frac{d}{p_0}$. Consider $\vu_0:\mathbb{R}^d\longrightarrow \mathbb{R}^d$ a divergence free initial data such that $\vu_0\in \dot{B}^{-(\alpha-1),\infty}_{\infty}(\mathbb{R}^d)$ and let $\vf:[0,+\infty[\times \mathbb{R}^d\longrightarrow \mathbb{R}^d$ be an external force such that
\begin{equation}\label{Def_NormeF1}
\|\vec{f}\|_{\mathcal{F}^{-\beta,p_0}_\rho}  = \underset{\tau>0}{\mathrm{sup}} \ \tau^\rho\|(-\Delta)^{\frac{-\beta}{2}}\vf(\tau, \cdot)\|_{L^{p_0}}<+\infty,\quad \mbox{with } \rho = 2-\tfrac{1}{\alpha}(\beta + \tfrac{d}{p_0}+1).
\end{equation}
In the case $p_0=+\infty$, we will ask moreover the extra condition $\mathrm{div}(\vec{f})=0$.\\ 
 
\noindent If the quantity $\|\vu_0\|_{\dot{B}^{-(\alpha-1),\infty}_{\infty}} + \|\vec{f}\|_{\mathcal{F}^{-\beta,p_0}_\rho}$ is small enough, then there exists a global \emph{mild} solution $\vu$ of the fractional integral problem (\ref{Formulation_Integrale}) such that $\vu \in L^{\infty}_\alpha([0,+\infty[\times \mathbb{R}^d)$ where 
$$L^{\infty}_\alpha([0,+\infty[\times\mathbb{R}^d)=\{\vpsi:[0,+\infty[\times\mathbb{R}^d\longrightarrow \mathbb{R}^d, \vpsi\in\mathcal{S'}([0,+\infty[\times \mathbb{R}^d):\|\vpsi\|_{L^{\infty}_\alpha}<+\infty\},$$ with 
\begin{equation}\label{Def_NormeX1}
\|\vpsi\|_{L^{\infty}_\alpha} = \underset{t>0}{\sup} \ t^{\frac{\alpha-1}{\alpha}}\|\vpsi(t, \cdot)\|_{L^{\infty}}.
\end{equation}
\end{Theoreme}
Some remarks are in order here. We first remark that the functional space $L^{\infty}_\alpha$ defined above is indeed critical in the sense that we have $\|\vu_\lambda\|_{L^{\infty}_\alpha}=\|\vu\|_{L^{\infty}_\alpha}$, which is the identity (\ref{Espaces_Critiques}) for the scaling given in (\ref{Homogeneite}). Note next that this space $L^{\infty}_\alpha$ was first considered in \cite{YuZhai} and the use of this particular space in the fixed-point argument allows to consider an initial data $\vu_0$ in the maximal critical Besov space $\dot{B}^{-(\alpha-1),\infty}_{\infty}$ (see Proposition \ref{Propo_EquivalenceBesov1} below). Remark also that we have the same criticality property for the space $\mathcal{F}^{-\beta,p_0}_\rho$ given by the norm (\ref{Def_NormeF1}) as we have $\|\vf_\lambda\|_{\mathcal{F}^{-\beta,p_0}_\rho}=\|\vf\|_{\mathcal{F}^{-\beta,p_0}_\rho}$ which corresponds with the identity (\ref{Espaces_CritiquesForce}) for the scaling (\ref{HomogeneiteForce}). Remark now that the condition $\div(\vf)=0$ when $p_0=+\infty$ is essentially related to the properties of the Leray projector and integrability issues. Finally, and to the best of our knowledge, this theorem gives a first existence result of mild solutions for the forced fractional Navier-Stokes equation (\ref{Formulation_Integrale}).\\

Although the space $L^{\infty}_\alpha$ characterized in (\ref{Def_NormeX1}) allows to close a fixed-point argument, more general funcional spaces can be considered and following the work \cite{quad_non} we have our next result which relies in the use of parabolic Morrey spaces. Recall that for a locally integrable function $\vpsi:[0,+\infty[\times\mathbb{R}^d\longrightarrow \mathbb{R}^d$ we have that $\vpsi$ belongs to the parabolic Morrey space $\mathcal{M}_\alpha^{p,q}([0,+\infty[\times\mathbb{R}^d)$ with $1<\alpha<2$ and $1\leq p \leq q <+\infty$ if 
\begin{equation}\label{Definition_MorreyPara}
\|\vpsi\|_{\mathcal{M}_\alpha^{p,q}}=\underset{r>0}{\mathrm{sup}} \ \underset{(t,x)\in [0,+\infty[ \times \mathbb{R}^d}{\mathrm{sup}} \ \frac{1}{r^{(d+\alpha)(\frac{1}{p}-\frac{1}{q})}} \left(\displaystyle{\iint_{\{|t-s|^\frac{1}{\alpha}+|x-y|<r\}}}|\vpsi(s,y)|^pdyds \right)^\frac{1}{p} <+\infty.
\end{equation}
Next and for some parameter $\gamma>0$, we will consider the following Morrey-based Sobolev space $\dot{\mathcal{W}}^{-\gamma, p,q}([0,+\infty[\times\mathbb{R}^d)$ given by the condition 
\begin{equation}\label{Definition_MorreySobolev}
\|\vpsi\|_{\dot{\mathcal{W}}^{-\gamma, p,q}}=\|(-\Delta)^{-\frac{\gamma}{2}}\vpsi\|_{\mathcal{M}_\alpha^{p,q}}<+\infty,
\end{equation}
(note that the operator $(-\Delta)^{-\frac{\gamma}{2}}$ acts only in the space variable). With these functional spaces we can state now the following theorem: 
\begin{Theoreme}\label{Theorem2}    
For $1<\alpha<2$ a fixed regularity index let $\vu_0:\mathbb{R}^d\longrightarrow \mathbb{R}^d$ be a divergence free initial data such that $\vu_0\in \dot{B}^{-(\alpha-1),\infty}_{\infty}(\mathbb{R}^d)$.\\

\noindent Consider now $0<p_1<+\infty$ a real parameter such that $2<p_1<\frac{\alpha}{\alpha-1}$
and define a real index $\gamma>0$ such that $2\alpha-1-(\alpha-1)p_1<\gamma<\alpha$.\\

\noindent Let $\vf:[0,+\infty[\times\mathbb{R}^d\longrightarrow \mathbb{R}^d$ be an external force such that we have $\vf \in \dot{\mathcal{W}}^{-\gamma, \mathfrak{p}, \mathfrak{q}}([0,+\infty[\times\mathbb{R}^d)$ with $\mathfrak{p}=\frac{(\alpha-1)p_1}{2\alpha-1-\gamma}>1$ and $\mathfrak{q}=\frac{d+\alpha}{2\alpha-1-\gamma}$, where this Morrey-based Sobolev space $\dot{\mathcal{W}}^{-\gamma, \mathfrak{p}, \mathfrak{q}}$ is characterized by the expression (\ref{Definition_MorreySobolev}). \\

\noindent If the quantity $\|\vu_0\|_{\dot{B}^{-(\alpha-1),\infty}_{\infty}} +\|\vec{f}\|_{\dot{\mathcal{W}}^{-\gamma, \mathfrak{p}, \mathfrak{q}}}$ is small enough, then there exists a global \emph{mild} solution $\vu$ of the fractional integral problem (\ref{Formulation_Integrale}) such that $\vu \in \mathcal{M}^{p_1,\frac{d+\alpha}{\alpha-1}}_\alpha([0,+\infty[\times \mathbb{R}^d)$ where 
 this parabolic Morrey space is characterized by the condition (\ref{Definition_MorreyPara}) above.
\end{Theoreme}
Let us make some observations. Remark first that, besides the parameter $1<\alpha<2$ and the dimension $d$ which are ``naturally'' given, the indexes $p_1, \gamma, \mathfrak{p}$ and $\mathfrak{q}$ that characterize the spaces $\dot{\mathcal{W}}^{-\gamma, \mathfrak{p}, \mathfrak{q}}$ and $\mathcal{M}^{p_1,\frac{d+\alpha}{\alpha-1}}_\alpha$ respond to different criteria. Indeed, for the resolution space $\mathcal{M}^{p_1,\frac{d+\alpha}{\alpha-1}}_\alpha$, the lower condition $2\leq p_1$ is related to some technical constraints (essentially it is a consequence of the H\"older inequalities in Morrey spaces, see also the expression (\ref{HolderMorreyPointFixe}) below) while the upper condition $p_1<\frac{\alpha}{\alpha-1}$ is stated to ensure the space inclusion $L^\infty_\alpha\subset \mathcal{M}^{p_1,\frac{d+\alpha}{\alpha-1}}_\alpha$ as well as the estimate $\|\mathfrak{p}_t\ast \vu_0\|_{ \mathcal{M}^{p_1,\frac{d+\alpha}{\alpha-1}}_\alpha}\leq C \|\vu_0\|_{\dot{B}^{-(\alpha-1),\infty}_\infty}$ (see more details in the Proposition \ref{Proposition_InclusiondansMorrey} and in the expression (\ref{EstimationDonneeInitiale1}) below). Note that the second index $\frac{d+\alpha}{\alpha-1}$ of the space $\mathcal{M}^{p_1,\frac{d+\alpha}{\alpha-1}}_\alpha$ guarantees the criticality of this resolution space as we have the homogeneity property $\|\vu_\lambda\|_{ \mathcal{M}^{p_1,\frac{d+\alpha}{\alpha-1}}_\alpha}=\|\vu\|_{ \mathcal{M}^{p_1,\frac{d+\alpha}{\alpha-1}}_\alpha}$ for the scaling given in (\ref{Homogeneite}). Remark also that for the space $\dot{\mathcal{W}}^{-\gamma, \mathfrak{p}, \mathfrak{q}}$, the indexes $\gamma$, $\mathfrak{p}$ and $\mathfrak{q}$ are given in order to obtain a control of the form (\ref{Controle_Force}) (see expression (\ref{EstimationMorreyForcePreuve}) below) and to obtain the identity $\|\vf_\lambda\|_{\dot{\mathcal{W}}^{-\gamma, \mathfrak{p}, \mathfrak{q}}}=\|\vf\|_{\dot{\mathcal{W}}^{-\gamma, \mathfrak{p}, \mathfrak{q}}}$ for the scaling given in (\ref{HomogeneiteForce}). We are thus working in a critical setting for the initial data, the external force and the resolution space. \\

Observe now that we have the following space inclusion (see Proposition \ref{Proposition_InclusiondansMorrey} below):
$$L^{\infty}_\alpha([0,+\infty[\times\mathbb{R}^d)\subset \mathcal{M}^{p_1,\frac{d+\alpha}{\alpha-1}}_\alpha([0,+\infty[\times\mathbb{R}^d),$$
and in this sense (at least when $\vf\equiv 0$) Theorem \ref{Theorem2} is a clear generalization of Theorem \ref{Theorem1}. However, the situation is slightly more subtle in the presence of an external force. Indeed, as it will be made explicit in Section \ref{Secc_ContreExemples} with some particular counterexamples, there is not a simple relationship between the spaces
$$\mathcal{F}^{-\beta,p_0}_\rho([0,+\infty[\times \mathbb{R}^d)\qquad \mbox{and}\qquad \dot{\mathcal{W}}^{-\gamma, \mathfrak{p}, \mathfrak{q}} ([0,+\infty[\times \mathbb{R}^d),$$
and this suggest that there it does not exist a ``simple'' criterion to consider the most general external force when studying the equation (\ref{EquationNS_intro}). Thus, in the presence of an external force, Theorem \ref{Theorem1} and Theorem \ref{Theorem2} can not be compared.\\[2mm]

We consider now another point of view. Indeed, by carefully studying the structure of the equation (\ref{Formulation_Integrale}) and by considering more suitable functional spaces, it is possible to ``generalize'' the previous results: we will work now with some particular multipliers spaces which are defined as follows: for $1<\alpha<2$, we first consider the subset $\mathcal{C}_\alpha$ of the local space $(L^2_tL^2_x)_{loc}$ defined by the condition
$$ v\in \mathcal{C}_\alpha \Longleftrightarrow  \int_{0}^{+\infty}\int_{\mathbb{R}^{d}} \frac{1}{(|t-s|^\frac{1}{\alpha}+|x-y|)^{d+1}}v(s,y)^2dyds\leq v(t,x),$$ 
where $v:[0,+\infty[\times \mathbb{R}^d\longrightarrow \mathbb{R}$ is a positive locally square integrable function. 
We will then say that a vector field $\vu:[0,+\infty[\times \mathbb{R}^d\longrightarrow \mathbb{R}^d$ belongs to the multiplier space $\mathcal{V}_\alpha([0,+\infty[\times \mathbb{R}^d, \mathbb{R}^d)$ if and only if there exist a constant $\mathfrak{C}\geq 0$ and a function $v\in \mathcal{C}_\alpha$ such that we have the (\emph{a.e.}) pointwise control:
\begin{equation}\label{Definition_Multiplicateurs1}
|\vec{u}(t,x)| \leq \mathfrak{C}\  v(t,x),
\end{equation}
we will thus define the quantity $\|\vu\|_{\mathcal{V}_\alpha}$ as the least constant $\mathfrak{C}=\mathfrak{C}(\vec{u})\geq 0$ in the previous estimate. In Section \ref{Secc_Notations} below we will present others characterization of the multiplier space $\mathcal{V}_\alpha$\\

With this functional space $\mathcal{V}_\alpha$, we can state our last result: 
\begin{Theoreme}\label{Theorem3}    
For $1<\alpha<2$ a fixed regularity index let $\vu_0:\mathbb{R}^d\longrightarrow \mathbb{R}^d$ be a divergence free initial data such that $\vu_0\in \dot{B}^{-(\alpha-1),\infty}_{\infty}(\mathbb{R}^d)$. Consider $\vf:[0,+\infty[\times \mathbb{R}^d\longrightarrow \mathbb{R}^d$ a given external force that belongs to the space $\mathcal{V}_\alpha^{-1}([0,+\infty[\times \mathbb{R}^d)$ which is defined by the condition
\begin{equation}\label{EspaceMultiplicateurForce0}
\vf \in \mathcal{V}_\alpha^{-1}([0,+\infty[\times \mathbb{R}^d)\iff\|\vec{f}\|_{\mathcal{V}_\alpha^{-1}}=\left\|\sqrt{|(-\Delta)^{-\frac{1}{2}}\vec{f}|}\right\|_{\mathcal{V}_\alpha}^2<+\infty.
\end{equation}
\noindent If the quantity $\|\vu_0\|_{\dot{B}^{-(\alpha-1),\infty}_{\infty}} +\|\vec{f}\|_{\mathcal{V}_\alpha^{-1}}$ is small enough, then there exists a global \emph{mild} solution $\vu$ of the fractional integral problem (\ref{Formulation_Integrale}) such that $\vu \in \mathcal{V}_\alpha([0,+\infty[\times \mathbb{R}^d)$. 
\end{Theoreme}
Note again that the space $\mathcal{V}_\alpha$ is critical in the sense that we have $\|\vu_\lambda\|_{\mathcal{V}_\alpha}=\|\vu\|_{\mathcal{V}_\alpha}$ for the scaling given in (\ref{Homogeneite}) and that the space $\mathcal{V}_\alpha^{-1}$ is also critical for the force since we have the identity $\|\vf_\lambda\|_{\mathcal{V}_\alpha^{-1}}=\|\vf\|_{\mathcal{V}_\alpha^{-1}}$ for the scaling (\ref{HomogeneiteForce}) and again we are working in a critical setting for the initial data, the external force and the resolution space. Remark now that, from the point of view of the resolution spaces, we have the following inclusion (see \cite[Theorem 4]{quad_non})
$$ \mathcal{M}^{p_1,\frac{d+\alpha}{\alpha-1}}_\alpha([0,+\infty[\times\mathbb{R}^d)\subset \mathcal{V}_\alpha([0,+\infty[\times\mathbb{R}^d),$$
and, at least when no forces are considered (\emph{i.e.} if $\vf\equiv 0$), this last result generalizes the Theorems \ref{Theorem1} \& \ref{Theorem2}. 
However, the situation is different when we take into account external forces. Indeed, as shown in Section \ref{Secc_ContreExemples}, below we do not have in general the inclusion of spaces
$$ \dot{\mathcal{W}}^{-\gamma, \mathfrak{p},\mathfrak{q}} ([0,+\infty[\times \mathbb{R}^d)\subset \mathcal{V}_\alpha^{-1}([0,+\infty[\times \mathbb{R}^d),$$
and thus, in presence of an external force, we can not compare, in general, Theorem \ref{Theorem2} with Theorem \ref{Theorem3}. \\

As we can observe, Theorem \ref{Theorem1}, Theorem \ref{Theorem2} and Theorem \ref{Theorem3} share some common features. First, all of them allow us to consider an initial data $\vu_0$ in the maximal Besov space $\dot{B}^{-(\alpha-1),\infty}_\infty$, which is in sharp contrast with the situation of the classical Navier-Stokes equations (\emph{i.e.} when $\alpha=2$) since the corresponding Besov space $\dot{B}^{-1,\infty}_\infty$ does not allow to close the fixed-point argument (recall the article \cite{besov_ill_posed}). Next, without external forces, there is a clear gain of generality as we have the inclusions   
$$L^{\infty}_\alpha([0,+\infty[\times\mathbb{R}^d)\subset  \mathcal{M}^{p_1,\frac{d+\alpha}{\alpha-1}}_\alpha([0,+\infty[\times\mathbb{R}^d)\subset \mathcal{V}_\alpha([0,+\infty[\times\mathbb{R}^d),$$
and thus the functional space $\mathcal{V}_\alpha$ seems to be -to the best of our knowledge- the largest space where to find \emph{mild} solutions for the problem (\ref{Formulation_Integrale}) from an initial data $\vu_0\in \dot{B}^{-(\alpha-1),\infty}_\infty$. The situation is however completely different in the presence of external forces. Indeed, there are not natural relationships between the spaces 
$$\mathcal{F}^{-\beta,p_0}_\rho([0,+\infty[\times \mathbb{R}^d),\qquad \dot{\mathcal{W}}^{-\gamma, \mathfrak{p}, \mathfrak{q}} ([0,+\infty[\times \mathbb{R}^d)\qquad \mbox{and}\qquad \mathcal{V}_\alpha^{-1}([0,+\infty[\times \mathbb{R}^d),$$
for \emph{all} the set of admissible parameters $\beta, p_0, \rho$, $\gamma, \mathfrak{p}, \mathfrak{q}$ that define these functional spaces and this shows the importance of each one of the previous theorems.\\

The plan of the article is the following. In Section \ref{Secc_Notations} we present some properties of the fractional heat kernel $\mathfrak{p}_t$ which are essential to our purposes and we recall the definition as well as some properties of the functional spaces used in the previous results. Sections \ref{Sec_Th1}, \ref{Sec_Th2} and \ref{Sec_Th3} are devoted to the proofs of the theorems stated above. In Section \ref{Secc_ContreExemples} we gave some counterexamples of the spaces used in the previous results for the external forces. 
\section{Fractional heat kernel and functional spaces}\label{Secc_Notations}
In this section we first gather some material related to the properties of the fractional heat kernel  $\mathfrak{p}_t$, with $1<\alpha<2$ and $t>0$, that will be crucial to perform our computations. Recall that we have $e^{-t(-\Delta)^{\frac{\alpha}{2}}}\vpsi=\mathfrak{p}_t\ast \vpsi$ (component wise) and thus, in the Fourier level, we have the identity $\widehat{\mathfrak{p}_t}(\xi)=e^{-t|\xi|^\alpha}$. For $1< \alpha < 2$, $\mathfrak{p}_t$ is a positive kernel and we have, for $C_1, C_2>0$:
$$C_1\frac{t}{(t^\frac{1}{\alpha}+|x|)^{d+\alpha}} \leq  \mathfrak{p}_t(x) \leq C_2\frac{t}{(t^\frac{1}{\alpha}+|x|)^{d+\alpha}},$$
a proof of this fact can be found in \cite[Theorem 7.3.1, p. 320]{kolokoltsov}. An easy consequence of these estimates is the inequality, valid for $1\leq p\leq +\infty$.
$$\|\mathfrak{p}_t\|_{L^{p}} \leq  Ct^{\frac{d}{\alpha p}-\frac{d}{\alpha}}.$$
Another crucial pointwise estimates are the ones obtained when a homogeneous pseudo-differential operator acts on the fractional heat kernel $\mathfrak{p}$. To this end, we introduce the following notation: let $\sigma \in \mathcal{C}^{\infty}(\mathbb{R}^d \setminus \{0\})$ be a positive, homogeneous function on $\mathbb{R}^d\setminus \{0\}$, with homogeneous degree $\gamma>-d$ (\emph{i.e.} we have $\sigma(\lambda \cdot)=\lambda^{\gamma}\sigma(\cdot)$ for all $\lambda>0$). 
For $1<\alpha<2$ and for $t>0$, we now define the function $\mathcal{K}^{\alpha,\sigma}_t:\mathbb{R}^d\longrightarrow \mathbb{R}$ in the Fourier level by the condition 
$$\widehat{\mathcal{K}^{\alpha,\sigma}_t}(\xi)=\sigma(\xi)e^{-t|\xi|^\alpha},$$
we will write moreover $\mathcal{K}^{\alpha,\sigma}$ when $\mathcal{K}^{\alpha,\sigma}_1$. It is worth mentioning that the condition $\gamma>-d$ is only used to ensure that $\sigma(\xi)e^{-t|\xi|^\alpha}$ is integrable with respect to $\xi$ in a neighborhood of the origin. We have the following homogeneity properties that are an easy consequences of the formula above: 
$$\widehat{\mathcal{K}^{\alpha,\sigma}_t}(\xi)=t^{-\frac{\gamma}{\alpha}} \widehat{\mathcal{K}^{\alpha,\sigma}}(t^\frac{1}{\alpha}\xi)\qquad \mbox{and}\qquad \mathcal{K}^{\alpha,\sigma}_t(x)= t^{-\frac{d+\gamma}{\alpha}} \mathcal{K}^{\alpha,\sigma} (t^\frac{-1}{\alpha}x).$$
We present now a pointwise estimate of such kernels, which will immediately grant us useful integrability properties: indeed, for $1<\alpha<2$, for $t>0$ and for some constant $C>0$, we have:
$$|\mathcal{K}^{\alpha,\sigma}_t(x)| \leq \frac{C}{(t^\frac{1}{\alpha}+|x|)^{d+\gamma}}.$$ 
In particular for $1\leq p \leq +\infty$, we obtain that $\mathcal{K}^{\alpha,\sigma}_t \in L^{p}(\mathbb{R}^d)$ whenever $\gamma > -\frac{d}{p'}$. Note that due to homogeneity properties we have
\begin{equation}\label{Estimation_NoyauFracc}
\left\|\mathcal{K}^{\alpha,\sigma}_t\right\|_{L^{p}} = t^{\frac{d}{\alpha p}-\frac{d+\gamma}{\alpha}} \left\|\mathcal{K}^{\alpha,\sigma}\right\|_{L^{p}}.
\end{equation}
See a proof of these facts in \cite{quad_non}. 
\subsubsection*{Functional spaces}
We will now fix here some notations and properties of the functional spaces used in this article. 
\begin{itemize}
\item[$\bullet$] We start with the Besov space for the initial data which will be used the most here and we will mainly use the thermic characterization for these spaces: indeed, we will say that $\psi\in \dot{B}^{-s,\infty}_\infty(\mathbb{R}^d)$ if 
$$\|\psi\|_{\dot{B}^{-s,\infty}_\infty}=\underset{t>0}{\sup}\; t^{\frac{s}{2}}\|h_t\ast \psi\|_{L^\infty}<+\infty,$$
where $h_t$ is the usual heat kernel. These spaces can also be defined using the fractional heat kernel $\mathfrak{p}_t$, indeed we have the following result
\begin{Proposition}\label{Propo_EquivalenceBesov1}
The Besov space $\dot{B}^{-s,\infty}_\infty(\mathbb{R}^d)$ can also be characterized by the condition 
$$\psi\in \dot{B}^{-s,\infty}_\infty(\mathbb{R}^d)\iff \underset{t>0}{\sup}\; t^{\frac{s}{\alpha}}\|\mathfrak{p}_t\ast \psi\|_{L^\infty}<+\infty,$$
where $\mathfrak{p}_t$ is the fractional heat kernel given by the expression (\ref{Definition_FracHeat}) above with $1<\alpha<2$. We thus have the equivalence 
$$ \underset{t>0}{\sup}\; t^{\frac{s}{2}}\|h_t\ast \psi\|_{L^\infty}\simeq  \underset{t>0}{\sup}\; t^{\frac{s}{\alpha}}\|\mathfrak{p}_t\ast \psi\|_{L^\infty}.$$
\end{Proposition}
See a proof of this equivalence in \cite[Proposition 2.1]{Miao}, see also \cite{quad_non}.\\

\item[$\bullet$] We turn now our attention to the parabolic Morrey spaces given in the expression (\ref{Definition_MorreyPara}) above. First note that if we set $p=q$ in (\ref{Definition_MorreyPara}), we obtain the identification $\mathcal{M}_\alpha^{p,q}([0,+\infty[\times \mathbb{R}^d)=L^p([0,+\infty[\times \mathbb{R}^d)$, recall also that we have the H\"older inequality
\begin{equation}\label{Holder_Morrey}
\|\vf\cdot \vec{g}\|_{\mathcal{M}^{p,q}_\alpha}\leq \|\vf\|_{\mathcal{M}^{p_1,q_1}_\alpha}\|\vec{g}\|_{\mathcal{M}^{p_2,q_2}_\alpha},
\end{equation}
where $1<p\leq q<+\infty$ with $\frac{1}{p}=\frac{1}{p_1}+\frac{1}{p_2}$ and $\frac{1}{q}=\frac{1}{q_1}+\frac{1}{q_2}$ (see \cite[Theorem 2.3]{Rafeiro}). Next, remark that by extending by $0$ a function $\vpsi: [0,+\infty[\times\mathbb{R}^{d}\longrightarrow \mathbb{R}^d$ we can consider without any problem the space $\mathcal{M}_\alpha^{p,q}(\mathbb{R}\times \mathbb{R}^d)$ and this extension by $0$ will allow us to study some classical (parabolic) operators. Indeed, for a suitable function $\vpsi:\mathbb{R}\times \mathbb{R}^d\longrightarrow \mathbb{R}^d$ we can consider the parabolic Riesz potential $\mathcal{I}_\mathfrak{s}(\vpsi)$ with $0<\mathfrak{s}<d+\alpha$ by the expression
\begin{equation}\label{DefParaboliqueRiesz}
\mathcal{I}_\mathfrak{s}(\vpsi)(t,x)=\int_\mathbb{R} \int_{\mathbb{R}^d} \frac{1}{(|t-s|^\frac{1}{\alpha}+|x-y|)^{d+\alpha-\mathfrak{s}}}|\vpsi(s,y)|dyds.
\end{equation}
We recall the following estimate: for a parameter $p$ such that $1<\frac{\mathfrak{s}+\alpha-1}{\alpha-1}<p\leq \frac{d+\alpha}{\alpha-1}$, we have 
\begin{equation}\label{ContinuiteRiesz}
\|\mathcal{I}_\mathfrak{s}(\vpsi)\|_{\mathcal{M}^{p,\frac{d+\alpha}{\alpha-1}}_\alpha}  \leq C\|\vpsi\|_{\mathcal{M}^{\mathfrak{p},\mathfrak{q}}_\alpha},
\end{equation}
where $\mathfrak{p}=\frac{(\alpha-1)p}{\mathfrak{s}+\alpha-1}$ and $\mathfrak{q}=\frac{d+\alpha}{\mathfrak{s}+\alpha-1}$. See a proof of this fact in \cite[Corollary 5.1]{PGLR}.
   
\item[$\bullet$] We give here some details about the space $\mathcal{V}_\alpha$. Recall that if we define $\|\vpsi\|_{\mathcal{V}_\alpha}$ as the least constant $\mathfrak{C}$ in the estimate (\ref{Definition_Multiplicateurs1}), then the space $\mathcal{V}_\alpha$ inherits a Banach function space structure (see \cite[Proposition 1]{quad_non}), in particular if we have $|\vpsi|\leq |\vec{\phi}|$ \emph{a.e.} then we have 
\begin{equation}\label{MonotiniciteNormeMultiplicateur}
\|\vpsi\|_{\mathcal{V}_\alpha}\leq \|\vec{\phi}\|_{\mathcal{V}_\alpha}.
\end{equation}
An important particularity from this space is that we have $\vpsi \in \mathcal{V}_\alpha$ if and only if $\mathcal{I}_{\alpha-1}(|\vpsi|^2) \in \mathcal{V}_\alpha$ and the we have the equivalence
\begin{equation}\label{EquivalenceMultiplicateur}
 \|\vpsi\|_{\mathcal{V}_\alpha} \simeq \sqrt{\left\| \mathcal{I}_{\alpha-1}(|\vpsi|^2) \right\|_{\mathcal{V}_\alpha}},
 \end{equation}
 where we recall that $\mathcal{I}_{\alpha-1}$ is the parabolic Riesz potential given in  (\ref{DefParaboliqueRiesz}). Another important feature of the multiplier space $\mathcal{V}_\alpha$ is that we have the control 
\begin{equation}\label{ControlRieszMultiplicateur}
\|\mathcal{I}_{\alpha-1}(|\vpsi| |\vec{\phi}|)\|_{\mathcal{V}_\alpha}\leq C\|\vpsi\|_{\mathcal{V}_\alpha}\|\vec{\phi}\|_{\mathcal{V}_\alpha}.
\end{equation}
For a proof of these facts and more details on the space $\mathcal{V}_\alpha$, see \cite{quad_non}.

\end{itemize}
\section{Proof of the Theorem \ref{Theorem1}}\label{Sec_Th1}
Recall that we are working here in the resolution space $(L^{\infty}_\alpha, \|\cdot\|_{L^{\infty}_\alpha})$ where the norm $\|\cdot\|_{L^{\infty}_\alpha}$ is given in the expression (\ref{Def_NormeX1}) above. 
It was proven in \cite[Theorem 3.1]{YuZhai} that for $1 < \alpha < 2$, we have the estimates 
$$\|\mathfrak{p}_t\ast\vu_0\|_{L^{\infty}_\alpha}\leq\|\vu_0\|_{\dot{B}^{-(\alpha-1)}_{\infty,\infty}},$$
as well as 
$$\left\|\int_0^t \mathfrak{p}_{t-s}\ast\mathbb{P}\left(\mathrm{div}(\vec{u}\otimes \vec{u})\right)ds\right\|_{L^{\infty}_\alpha}\leq C\|\vu\|_{L^{\infty}_\alpha}\|\vu\|_{L^{\infty}_\alpha},$$
which are the controls (\ref{Controle_DonneeInitiale}) and (\ref{Controle_Bilineaire}). The estimate for the initial data $\vu_0$ is a consequence of the Proposition \ref{Propo_EquivalenceBesov1} which actually states the equivalence of these quantities. For the bilinear term, we briefly recall here the arguments used in \cite{YuZhai}. Indeed we have 
\begin{eqnarray*}
\left\|\int_0^t \mathfrak{p}_{t-s}\ast\mathbb{P}\left(\mathrm{div}(\vec{u}\otimes \vec{u})\right)ds\right\|_{L^{\infty}_\alpha}&=&\underset{t>0}{\sup} \ t^{\frac{\alpha-1}{\alpha}}\left\|\int_0^t \mathfrak{p}_{t-s}\ast\mathbb{P}\left(\mathrm{div}(\vec{u}\otimes \vec{u})\right)ds\right\|_{L^\infty}\\
&\leq & \underset{t>0}{\sup} \ t^{\frac{\alpha-1}{\alpha}}\int_0^t\left\| \mathfrak{p}_{t-s}\ast\mathbb{P}\left(\mathrm{div}(\vec{u}\otimes \vec{u})\right)\right\|_{L^\infty}ds.
\end{eqnarray*}
At this point we use the following pointwise estimate\footnote{This pointwise estimate remains true if we replace the divergence operator by a generic derivative, say \emph{e.g.} $(-\Delta)^{\frac{1}{2}}$.} (see \cite[Lemma 1]{quad_non}) which is valid for a tensor $\mathbb{F}$ such that $\mathrm{div}(\mathbb{F}):\mathbb{R}^d\longrightarrow \mathbb{R}^d$:
\begin{equation}\label{fracOseen}
|\mathfrak{p}_{t-s}\ast\mathbb{P}(\mathrm{div}(\mathbb{F}))(x)|\leq C\int_{\mathbb{R}^d}\frac{|\mathbb{F}(y)|}{(|t-s|^{\frac{1}{\alpha}}+|x-y|)^{d+1}}dy.
\end{equation}
We can thus write
$$\left\|\int_0^t \mathfrak{p}_{t-s}\ast\mathbb{P}\left(\mathrm{div}(\vec{u}\otimes \vec{u})\right)ds\right\|_{L^{\infty}_\alpha}\leq C \underset{t>0}{\sup} \ t^{\frac{\alpha-1}{\alpha}}\int_0^t \left(\int_{\mathbb{R}^d}\frac{|\vu(s,y)||\vu(s,y)|}{(|t-s|^{\frac{1}{\alpha}}+|x-y|)^{d+1}}dy\right)ds,$$
from which we easily obtain the inequality 
\begin{eqnarray*}
\left\|\int_0^t \mathfrak{p}_{t-s}\ast\mathbb{P}\left(\mathrm{div}(\vec{u}\otimes \vec{u})\right)ds\right\|_{L^{\infty}_\alpha}&\leq &C \underset{t>0}{\sup} \ t^{\frac{\alpha-1}{\alpha}}\int_0^t \|\vu(s,\cdot)\|_{L^\infty}\|\vu(s,\cdot)\|_{L^\infty}\\
&&\times \left(\int_{\mathbb{R}^d}\frac{1}{(|t-s|^{\frac{1}{\alpha}}+|x-y|)^{d+1}}dy\right)ds\\[2mm]
&\leq & C \underset{t>0}{\sup} \ t^{\frac{\alpha-1}{\alpha}}\int_0^t \|\vu(s,\cdot)\|_{L^\infty}\|\vu(s,\cdot)\|_{L^\infty}(t-s)^{-\frac{1}{\alpha}}ds.
\end{eqnarray*}
Now, introducing the weight $s^{\frac{\alpha-1}{\alpha}}$ in the previous integral we obtain
\begin{eqnarray*}
\left\|\int_0^t \mathfrak{p}_{t-s}\ast\mathbb{P}\left(\mathrm{div}(\vec{u}\otimes \vec{u})\right)ds\right\|_{L^{\infty}_\alpha}&\leq &C \underset{t>0}{\sup} \ t^{\frac{\alpha-1}{\alpha}}\int_0^t (t-s)^{-\frac{1}{\alpha}}s^{-2(\frac{\alpha-1}{\alpha})} \\
&& \times (s^{\frac{\alpha-1}{\alpha}}\|\vu(s, \cdot)\|_{L^\infty})(s^{\frac{\alpha-1}{\alpha}}\|\vu(s, \cdot)\|_{L^\infty})ds\\[2mm]
&\leq & C\|\vu\|_{L^{\infty}_\alpha}\|\vu\|_{L^{\infty}_\alpha}\quad \underset{t>0}{\sup} \ t^{\frac{\alpha-1}{\alpha}}\int_0^t (t-s)^{-\frac{1}{\alpha}}s^{-2(\frac{\alpha-1}{\alpha})} ds,
\end{eqnarray*}
where we used the definition of the quantity $\|\cdot\|_{L^{\infty}_\alpha}$ given in (\ref{Def_NormeX1}). Now, a straightforward computation gives $\displaystyle{\int_0^t (t-s)^{-\frac{1}{\alpha}}s^{-2(\frac{\alpha-1}{\alpha})} ds}=Ct^{\frac{1-\alpha}{\alpha}}$, so we finally obtain 
$$\left\|\int_0^t \mathfrak{p}_{t-s}\ast\mathbb{P}\left(\mathrm{div}(\vec{u}\otimes \vec{u})\right)ds\right\|_{L^{\infty}_\alpha}\leq C\|\vu\|_{L^{\infty}_\alpha}\|\vu\|_{L^{\infty}_\alpha},$$
which corresponds with the boundedness (\ref{Controle_Bilineaire}) of the bilinear application in the space $L^\infty_\alpha$.\\[5mm]

To end the proof of Theorem \ref{Theorem1}, we only need to prove now the inequality (\ref{Controle_Force}), \emph{i.e.}:
\begin{equation}\label{Estimation_Force1}
\left\|\int_0^t \mathfrak{p}_{t-s}\ast\mathbb{P}(\vec{f})ds\right\|_{L^{\infty}_\alpha}= \underset{t>0}{\sup} \ t^{\frac{\alpha-1}{\alpha}}\left\|\int_0^t \mathfrak{p}_{t-s}\ast\mathbb{P}(\vf)ds\right\|_{L^{\infty}}\leq C\|\vf\|_{\mathcal{F}^{-\beta,p_0}_\rho},
\end{equation}
where $\|\vec{f}\|_{\mathcal{F}^{-\beta,p_0}_\rho}  = \underset{\tau>0}{\mathrm{sup}} \ \tau^\rho\|(-\Delta)^{-\frac{\beta}{2}}\vf(\tau, \cdot)\|_{L^{p_0}}$. Here the indexes $\sigma, \gamma$ and $p_0$ are related to the dimension $d$ and to the regularity parameter $1<\alpha<2$ and satisfy the conditions $\frac{d}{\alpha}<p_0 \leq +\infty$,  $\alpha-\frac{d}{p_0}-1< \beta < \alpha-\frac{d}{p_0}$ and $\rho = 2-\tfrac{1}{\alpha}(\beta + \tfrac{d}{p_0}+1)$, from which we can see that the choice $d\geq2$ and $\alpha<2$ guarantees that $p_0>\frac{d}{\alpha}>1$ and this implies that the interval in which $\beta$ is chosen is non empty.\\

In order to establish the estimate (\ref{Estimation_Force1}) we will decompose our analysis following the values of the parameter $p_0$.  
\begin{itemize}
\item Let us first consider the case $\frac{d}{\alpha}<p_0<+\infty$. We thus write
$$\left\|\int_0^t \mathfrak{p}_{t-s}\ast\mathbb{P}(\vf)(s, \cdot)ds\right\|_{L^{\infty}}\leq \int_0^t \left\|\mathfrak{p}_{t-s}\ast\mathbb{P}(\vf)(s, \cdot)\right\|_{L^{\infty}}ds,$$
and by the Young inequalities for convolution, the properties of the Leray projector, as well as by the properties of the fractional Laplacian, we have
\begin{eqnarray*}
\left\|\int_0^t \mathfrak{p}_{t-s}\ast\mathbb{P}(\vf)(s, \cdot)ds\right\|_{L^{\infty}}&\leq &\int_0^t \|(-\Delta)^{\frac{\beta}{2}}\mathfrak{p}_{t-s}\|_{L^{p_0'}}\|\mathbb{P}((-\Delta)^{\frac{-\beta}{2}}(\vf)(s, \cdot))\|_{L^{p_0}}ds\\
&\leq &C\int_0^t \|(-\Delta)^{\frac{\beta}{2}}\mathfrak{p}_{t-s}\|_{L^{p_0'}}\|(-\Delta)^{\frac{-\beta}{2}}(\vf)(s, \cdot)\|_{L^{p_0}}ds,
\end{eqnarray*}
where in the last estimate above we used the boundedness of the Leray projector in Lebesgue spaces $L^{p_0}$ with $1<p_0<+\infty$. Now by the properties of the fractional heat kernel given in (\ref{Estimation_NoyauFracc}) we obtain
\begin{eqnarray*}
\left\|\int_0^t \mathfrak{p}_{t-s}\ast\mathbb{P}(\vf)(s, \cdot)ds\right\|_{L^{\infty}}&\leq &C\int_0^t(t-s)^{-\frac{d}{\alpha p_0}-\frac{\beta}{\alpha}}\quad s^{-\rho} s^\rho \|(-\Delta)^{\frac{-\beta}{2}}(\vf)(s, \cdot)\|_{L^{p_0}}ds\\
&\leq &C\int_0^t(t-s)^{-\frac{d}{\alpha p_0}-\frac{\beta}{\alpha}} \;s^{-\rho} ds\quad \left(\underset{s>0}{\sup}\, s^\rho \|(-\Delta)^{\frac{-\beta}{2}}(\vf)(s, \cdot)\|_{L^{p_0}}\right)\\
&\leq &C\|\vf\|_{\mathcal{F}^{-\beta,p_0}_\rho}\int_0^t(t-s)^{-\frac{d}{\alpha p_0}-\frac{\beta}{\alpha}} s^{-\rho} ds,
\end{eqnarray*}
where in the lines above we introduced the quantity $\|\cdot\|_{\mathcal{F}^{-\beta,p_0}_\rho}$ given in (\ref{Def_NormeF1}). We write now 
$$\int_0^t(t-s)^{-\frac{d}{\alpha p_0}-\frac{\beta}{\alpha}} s^{-\rho} ds=\int_0^{\frac{t}{2}}(t-s)^{-\frac{d}{\alpha p_0}-\frac{\beta}{\alpha}} s^{-\rho} ds+\int_{\frac{t}{2}}^t(t-s)^{-\frac{d}{\alpha p_0}-\frac{\beta}{\alpha}} s^{-\rho} ds,$$
for the first integral we note that if $0\leq s\leq \frac{t}{2}$, then we have $\frac{t}{2}\leq (t-s)\leq t$ from which we obtain $(t-s)^{-\frac{d}{\alpha p_0}-\frac{\beta}{\alpha}}\leq Ct^{-\frac{d}{\alpha p_0}-\frac{\beta}{\alpha}}$ and if $\frac{t}{2}<s\leq t$ we have $s^{-\rho}\leq Ct^{-\rho}$, so we can write 
$$\int_0^t(t-s)^{-\frac{d}{\alpha p_0}-\frac{\beta}{\alpha}} s^{-\rho} ds\leq Ct^{-\frac{d}{\alpha p_0}-\frac{\beta}{\alpha}} \int_0^{\frac{t}{2}}s^{-\rho} ds+Ct^{-\rho}\int_{\frac{t}{2}}^t(t-s)^{-\frac{d}{\alpha p_0}-\frac{\beta}{\alpha}} ds,$$
and after an integration (recall that $\rho=2-\tfrac{1}{\alpha}(\beta + \tfrac{d}{p_0}+1)<1$ as $\alpha-\frac{d}{p_0}-1<\beta$ and $\frac{d}{\alpha p_0}+\frac{\beta}{\alpha}<1$ as $\beta<\alpha-\frac{d}{p_0}$) we thus have $\displaystyle{\int_0^t(t-s)^{-\frac{d}{\alpha p_0}-\frac{\beta}{\alpha}} s^{-\rho} ds\leq Ct^{1-\rho-\frac{d}{\alpha p_0}-\frac{\beta}{\alpha}}}$. 
Using the fact that $\rho = 2-\tfrac{1}{\alpha}(\beta + \tfrac{d}{p_0}+1)$ we obtain 
$\displaystyle{\int_0^t(t-s)^{-\frac{d}{\alpha p_0}-\frac{\beta}{\alpha}} s^{-\rho} ds\leq Ct^{-\frac{\alpha-1}{\alpha}}}$ and we can write
$$\left\|\int_0^t \mathfrak{p}_{t-s}\ast\mathbb{P}(\vf)(s, \cdot)ds\right\|_{L^{\infty}}\leq Ct^{-\frac{\alpha-1}{\alpha}} \|\vf\|_{\mathcal{F}^{-\beta,p_0}_\rho},$$
from which we easily deduce the wished estimate (\ref{Estimation_Force1}). 

\item When $p_0=+\infty$ we are assuming that $\vec{f}$ is divergence free and since we have $\mathbb{P}(\vf)=\vf$, we can write 
$$\left\|\int_0^t \mathfrak{p}_{t-s}\ast\mathbb{P}(\vf)(s, \cdot)ds\right\|_{L^{\infty}}=\left\|\int_0^t \mathfrak{p}_{t-s}\ast\vf (s, \cdot)ds\right\|_{L^{\infty}}\leq \int_0^t \left\| \mathfrak{p}_{t-s}\ast\vf (s, \cdot)\right\|_{L^{\infty}}ds,$$
and we obtain 
$$\left\|\int_0^t \mathfrak{p}_{t-s}\ast\mathbb{P}(\vf)(s, \cdot)ds\right\|_{L^{\infty}}\leq \int_0^t\|(-\Delta)^{\frac{\beta}{2}}\mathfrak{p}_{t-s}\|_{L^1}\|(-\Delta)^{-\frac{\beta}{2}}\vf(s, \cdot)\|_{L^\infty}ds,$$
from which we deduce, by the smoothness properties of the fractional heat kernel given in (\ref{Estimation_NoyauFracc}): 
$$\left\|\int_0^t \mathfrak{p}_{t-s}\ast\mathbb{P}(\vf)(s, \cdot)ds\right\|_{L^{\infty}}\leq C\int_0^t (t-s)^{-\frac{d+\beta}{\alpha}}\|(-\Delta)^{-\frac{\beta}{2}}\vf(s, \cdot)\|_{L^\infty}ds.$$
We introduce now the quantity $\|\cdot\|_{\mathcal{F}^{-\beta,\infty}_\rho}$ given in (\ref{Def_NormeF1}) to obtain
\begin{eqnarray*}
\left\|\int_0^t \mathfrak{p}_{t-s}\ast\mathbb{P}(\vf)(s, \cdot)ds\right\|_{L^{\infty}}&\leq &C\int_0^t (t-s)^{-\frac{d+\beta}{\alpha}} s^{-\rho} ds \left(\underset{s>0}{\sup}\; s^{\rho}\|(-\Delta)^{-\frac{\beta}{2}}\vf(s, \cdot)\|_{L^\infty}\right)\\
&\leq & C \|\vf\|_{\mathcal{F}^{-\beta,\infty}_\rho}\int_0^t (t-s)^{-\frac{d+\beta}{\alpha}} s^{-\rho} ds.
\end{eqnarray*}
Recall now that in the case $p_0=+\infty$ we have  $\rho=2-\tfrac{1}{\alpha}(\beta +1)<1$ and by the same arguments as above we obtain $\displaystyle{\int_0^t (t-s)^{-\frac{d+\beta}{\alpha}} s^{-\rho} ds}\leq Ct^{-\frac{(\alpha-1)}{\alpha}}$, we thus have 
$$\left\|\int_0^t \mathfrak{p}_{t-s}\ast\mathbb{P}(\vf)(s, \cdot)ds\right\|_{L^{\infty}}\leq Ct^{-\frac{(\alpha-1)}{\alpha}} \|\vf\|_{\mathcal{F}^{-\beta,\infty}_\rho},$$
which is the estimate (\ref{Estimation_Force1}) in the case $p_0=+\infty$.
\end{itemize}
The proof of Theorem \ref{Theorem1} is now complete. \hfill$\blacksquare$
\section{Proof of the Theorem \ref{Theorem2}}\label{Sec_Th2}
We work now in the resolution space $\mathcal{M}_{\alpha}^{p_1, \frac{d+\alpha}{\alpha-1}}([0,+\infty[\times \mathbb{R}^d)$ which is determined by the quantity (\ref{Definition_MorreyPara}) given above. As before, we need to prove the estimates
\begin{equation}\label{EstimationMorreyDonneeInitiale}
\|\mathfrak{p}_t\ast\vu_0\|_{\mathcal{M}_{\alpha}^{p_1, \frac{d+\alpha}{\alpha-1}}}\leq C_2\|\vu_0\|_{\dot{B}^{-(\alpha-1)}_{\infty,\infty}},
\end{equation}
for the initial data $\vu$, and
\begin{equation}\label{EstimationMorreyBilineaire}
\left\|\int_0^t \mathfrak{p}_{t-s}\ast\mathbb{P}\left(\mathrm{div}(\vec{u}\otimes \vec{u})\right)ds\right\|_{\mathcal{M}_{\alpha}^{p_1, \frac{d+\alpha}{\alpha-1}}}\leq C_2\|\vu\|_{\mathcal{M}_{\alpha}^{p_1, \frac{d+\alpha}{\alpha-1}}}\|\vu\|_{\mathcal{M}_{\alpha}^{p_1, \frac{d+\alpha}{\alpha-1}}},
\end{equation}
for the bilinear term, and 
\begin{equation}\label{EstimationMorreyForce}
\left\|\int_0^t \mathfrak{p}_{t-s}\ast\mathbb{P}(\vf)ds\right\|_{\mathcal{M}_{\alpha}^{p_1, \frac{d+\alpha}{\alpha-1}}}\leq C_3\|\vf\|_{\dot{\mathcal{W}}_\alpha^{-\rho, \mathfrak{p}, \mathfrak{q}}},
\end{equation}
for the external force $\vf$. We will  study each one of these terms separately, but before proving these estimates, we will show in the next proposition that the functional space $\mathcal{M}_{\alpha}^{p_1, \frac{d+\alpha}{\alpha-1}}([0,+\infty[\times \mathbb{R}^d)$ is indeed bigger that the resolution space $L^{\infty}_\alpha([0,+\infty[\times \mathbb{R}^d)$ considered in the previous section (and in the article \cite{YuZhai}), and in this sense we can claim that Theorem \ref{Theorem2} is a generalization of the Theorem \ref{Theorem1}. 
\begin{Proposition}\label{Proposition_InclusiondansMorrey}
Let $1\leq p_1<\frac{\alpha}{\alpha-1}$, then we have the space inclusion 
$$L^{\infty}_\alpha([0,+\infty[\times \mathbb{R}^d)\subset\mathcal{M}_{\alpha}^{p_1, \frac{d+\alpha}{\alpha-1}}([0,+\infty[\times \mathbb{R}^d).$$
\end{Proposition}
{\bf Proof. } To prove this fact we will consider a function $\vpsi:[0,+\infty[\times \mathbb{R}^d\longrightarrow \mathbb{R}^d$ and we will prove the estimate
$$\|\vpsi\|_{\mathcal{M}_{\alpha}^{p_1, \frac{d+\alpha}{\alpha-1}}}\leq C\|\vpsi\|_{L^{\infty}_\alpha}.$$
To this end, we recall that we have (see formula (\ref{Definition_MorreyPara}) above):
$$\|\vpsi\|_{\mathcal{M}_\alpha^{p_1, \frac{d+\alpha}{\alpha-1}}}=\underset{r>0}{\mathrm{sup}} \ \underset{(t,x)\in [0,+\infty[ \times \mathbb{R}^d}{\mathrm{sup}} \ \frac{1}{r^{(d+\alpha)(\frac{1}{p_1}- \frac{\alpha-1}{d+\alpha})}} \left(\displaystyle{\iint_{\{|t-s|^\frac{1}{\alpha}+|x-y|<r\}}}|\vpsi(s,y)|^{p_1}dyds \right)^\frac{1}{p_1},$$
and we write 
\begin{eqnarray*}
\|\vpsi\|_{\mathcal{M}_\alpha^{p_1, \frac{d+\alpha}{\alpha-1}}}=\hspace{14cm}\\
\underset{r>0}{\mathrm{sup}} \ \underset{(t,x)\in [0,+\infty[ \times \mathbb{R}^d}{\mathrm{sup}} \ \frac{1}{r^{(d+\alpha)(\frac{1}{p_1}-\frac{\alpha-1}{d+\alpha})}} \left(\displaystyle{\iint_{\{|t-s|^\frac{1}{\alpha}+|x-y|<r\}}}|s|^{-p_1(\frac{\alpha-1}{\alpha})}|s|^{p_1(\frac{\alpha-1}{\alpha})}|\vpsi(s,y)|^{p_1}dyds \right)^\frac{1}{p_1}\\
\leq  \underset{s>0}{\sup}\; s^{\frac{\alpha-1}{\alpha}}\|\vpsi(s,\cdot)\|_{L^\infty} \times \underset{r>0}{\mathrm{sup}} \ \underset{(t,x)\in [0,+\infty[ \times \mathbb{R}^d}{\mathrm{sup}} \ \frac{1}{r^{(d+\alpha)(\frac{1}{p_1}-\frac{\alpha-1}{d+\alpha})}} \left(\displaystyle{\iint_{\{|t-s|^\frac{1}{\alpha}+|x-y|<r\}}}|s|^{-p_1(\frac{\alpha-1}{\alpha})}dyds \right)^\frac{1}{p_1}
\end{eqnarray*}
\begin{eqnarray} \label{Estimation_pourInclusionMorreyChinois}
\leq \|\vpsi\|_{L^\infty_\alpha}\times \underset{r>0}{\mathrm{sup}} \ \underset{(t,x)\in [0,+\infty[ \times \mathbb{R}^d}{\mathrm{sup}} \ \frac{1}{r^{(d+\alpha)(\frac{1}{p_1}-\frac{\alpha-1}{d+\alpha})}} \left(\displaystyle{\iint_{\{|t-s|^\frac{1}{\alpha}+|x-y|<r\}}}|s|^{-p_1(\frac{\alpha-1}{\alpha})}dyds \right)^\frac{1}{p_1},
\end{eqnarray}
where we used the fact that $\|\vpsi\|_{L^\infty_\alpha}= \underset{s>0}{\sup}\; s^{\frac{\alpha-1}{\alpha}}\|\vpsi(s,\cdot)\|_{L^\infty}$. To conclude, we only need to study the integrals above. Since the set $\{|t-s|^\frac{1}{\alpha}+|x-y|<r\}$ is included in the set   $\{s>0: |t-s|<r^\alpha\} \times \{y\in \mathbb{R}^d: |x-y|<r\}$ we have
$$I=\iint_{\{|t-s|^\frac{1}{\alpha}+|x-y|<r\}}|s|^{-p_1(\frac{\alpha-1}{\alpha})}dyds\leq C r^{d}\int_{\{|t-s|<r^\alpha\}}s^{-p_1(\frac{\alpha-1}{\alpha})}ds=Cr^{d}\int_{t-r^\alpha}^{t+r^\alpha}|s|^{-p_1(\frac{\alpha-1}{\alpha})}ds.$$
We decompose our study of the last integral in two cases. First, if $0\leq t\leq2r^\alpha$, the domain of integration is then included in the interval $[-r^\alpha, 3r^\alpha]$ and we can write
$$I=Cr^d\int_{-r^\alpha}^{3r^\alpha}|s|^{-p_1\frac{\alpha-1}{\alpha}}ds\leq Cr^d\int_{-3r^\alpha}^{3 r^\alpha}|s|^{-p_1\frac{\alpha-1}{\alpha}}ds=2Cr^d\int_{0}^{3 r^\alpha}s^{-p_1\frac{\alpha-1}{\alpha}}ds,$$
since $1\leq p_1<\frac{\alpha}{\alpha-1}$, the previous integral is finite and we have 
$$I\leq C r^d r^{\alpha(1-p_1\frac{\alpha-1}{\alpha})}=Cr^{(d+\alpha)(1-p_1\frac{\alpha-1}{d+\alpha})}.$$
We consider now the case when $t>2r^\alpha$. Since we have here $t+r^\alpha>t-r^\alpha>r^\alpha>0$ and $-p_1\frac{\alpha-1}{\alpha}<0$, we can write
$$I\leq Cr^{d}\int_{t-r^\alpha}^{t+r^\alpha}|s|^{-p_1(\frac{\alpha-1}{\alpha})}ds\leq Cr^{d}\int_{t-r^\alpha}^{t+r^\alpha}(r^\alpha)^{-p_1(\frac{\alpha-1}{\alpha})}ds=Cr^{d-p_1(\alpha-1)}r^\alpha=Cr^{(d+\alpha)(1-p_1\frac{\alpha-1}{d+\alpha})}.$$
With these estimates for $I$ we come back to the formula (\ref{Estimation_pourInclusionMorreyChinois}) and we obtain 
\begin{eqnarray*}
\|\vpsi\|_{\mathcal{M}_\alpha^{p_1, \frac{d+\alpha}{\alpha-1}}}&\leq &\|\vpsi\|_{L^\infty_\alpha}\times \underset{r>0}{\mathrm{sup}} \ \underset{(t,x)\in [0,+\infty[ \times \mathbb{R}^d}{\mathrm{sup}} \ \frac{1}{r^{(d+\alpha)(\frac{1}{p_1}-\frac{\alpha-1}{d+\alpha})}} \left(\displaystyle{\iint_{\{|t-s|^\frac{1}{\alpha}+|x-y|<r\}}}|s|^{-p_1(\frac{\alpha-1}{\alpha})}dyds \right)^\frac{1}{p_1}\\
&\leq & C \|\vpsi\|_{L^\infty_\alpha}\times\underset{r>0}{\mathrm{sup}} \ \underset{(t,x)\in [0,+\infty[ \times \mathbb{R}^d}{\mathrm{sup}} \ \frac{1}{r^{(d+\alpha)(\frac{1}{p_1}-\frac{\alpha-1}{d+\alpha})}} r^{(d+\alpha)(\frac{1}{p_1}-\frac{\alpha-1}{d+\alpha})}\\
&\leq & C\|\vpsi\|_{L^\infty_\alpha},
\end{eqnarray*}
which is the wished inequality and this ends the proof of the Proposition \ref{Proposition_InclusiondansMorrey}. \hfill$\blacksquare$
\begin{Remarque}
Note that the condition $1\leq p_1<\frac{\alpha}{\alpha-1}$ is used here to ensure that the integral $I$ given in the previous computations converges.
\end{Remarque}
Once we have established that the functional space $\mathcal{M}_{\alpha}^{p_1, \frac{d+\alpha}{\alpha-1}}([0,+\infty[\times \mathbb{R}^d)$ is indeed bigger that the space $L^{\infty}_\alpha([0,+\infty[\times \mathbb{R}^d)$, we can now study the inequalities (\ref{EstimationMorreyDonneeInitiale}), (\ref{EstimationMorreyBilineaire}) and (\ref{EstimationMorreyForce}).

\begin{itemize}
\item For the initial data, \emph{i.e.} for the inequality (\ref{EstimationMorreyDonneeInitiale}),  we have: 
$$\|\mathfrak{p}_t\ast\vu_0\|_{\mathcal{M}_{\alpha}^{p_1, \frac{d+\alpha}{\alpha-1}}}=\underset{r>0}{\mathrm{sup}} \ \underset{(t,x)\in [0,+\infty[ \times \mathbb{R}^d}{\mathrm{sup}} \ \frac{1}{r^{(d+\alpha)(\frac{1}{p_1}- \frac{\alpha-1}{d+\alpha})}} \left(\displaystyle{\iint_{\{|t-s|^\frac{1}{\alpha}+|x-y|<r\}}}|\mathfrak{p}_s\ast\vu_0(y)|^{p_1}dyds \right)^\frac{1}{p_1},$$
and we can write
\begin{eqnarray*}
\|\mathfrak{p}_t\ast\vu_0\|_{\mathcal{M}_{\alpha}^{p_1, \frac{d+\alpha}{\alpha-1}}}&\leq & \underset{s>0}{\sup}\; s^{\frac{\alpha-1}{\alpha}}\|\mathfrak{p}_s\ast\vu_0\|_{L^\infty} \times \\ 
&&\underset{r>0}{\mathrm{sup}} \ \underset{(t,x)\in [0,+\infty[ \times \mathbb{R}^d}{\mathrm{sup}} \ \frac{1}{r^{(d+\alpha)(\frac{1}{p_1}- \frac{\alpha-1}{d+\alpha})}} \left(\displaystyle{\iint_{\{|t-s|^\frac{1}{\alpha}+|x-y|<r\}}}|s|^{-p_1(\frac{\alpha-1}{\alpha})} dyds \right)^\frac{1}{p_1}.
\end{eqnarray*}
Using Proposition \ref{Propo_EquivalenceBesov1} above we have $\underset{s>0}{\sup}\; s^{\frac{\alpha-1}{\alpha}}\|\mathfrak{p}_s\ast\vu_0\|_{L^\infty}\simeq \|\vu_0\|_{\dot{B}^{-(\alpha-1),\infty}_\infty}$ and by the computations performed in the previous proposition (see formula (\ref{Estimation_pourInclusionMorreyChinois})), we easily obtain the control 
\begin{equation}\label{EstimationDonneeInitiale1}
\|\mathfrak{p}_t\ast\vu_0\|_{\mathcal{M}_{\alpha}^{p_1, \frac{d+\alpha}{\alpha-1}}}\leq C\|\vu_0\|_{\dot{B}^{-(\alpha-1),\infty}_\infty}.
\end{equation}
\item For the bilinear term estimate (\ref{EstimationMorreyBilineaire}) we write:
\begin{eqnarray*}
\left\|\int_0^{\tau} \mathfrak{p}_{\tau-s}\ast\mathbb{P}\left(\mathrm{div}(\vec{u}\otimes \vec{u})\right)ds\right\|_{\mathcal{M}_{\alpha}^{p_1, \frac{d+\alpha}{\alpha-1}}}=\underset{r>0}{\mathrm{sup}} \ \underset{(t,x)\in [0,+\infty[ \times \mathbb{R}^d}{\mathrm{sup}} \ \frac{1}{r^{(d+\alpha)(\frac{1}{p_1}- \frac{\alpha-1}{d+\alpha})}}\qquad \qquad\\
\times \left(\displaystyle{\iint_{\{|t-\tau|^\frac{1}{\alpha}+|x-y|<r\}}}\left| \int_0^{\tau} \mathfrak{p}_{\tau-s}\ast\mathbb{P}\left(\mathrm{div}(\vec{u}\otimes \vec{u})\right)ds\right|^{p_1}dyd\tau \right)^\frac{1}{p_1},
\end{eqnarray*}
and using the pointwise estimate (\ref{fracOseen}) we have
\begin{eqnarray*}
\left\|\int_0^\tau \mathfrak{p}_{\tau-s}\ast\mathbb{P}\left(\mathrm{div}(\vec{u}\otimes \vec{u})\right)ds\right\|_{\mathcal{M}_{\alpha}^{p_1, \frac{d+\alpha}{\alpha-1}}}\leq C\underset{r>0}{\mathrm{sup}} \ \underset{(t,x)\in [0,+\infty[ \times \mathbb{R}^d}{\mathrm{sup}} \ \frac{1}{r^{(d+\alpha)(\frac{1}{p_1}- \frac{\alpha-1}{d+\alpha})}}\qquad \qquad\\
\times \left(\displaystyle{\iint_{\{|t-\tau|^\frac{1}{\alpha}+|x-y|<r\}}}\left( \int_0^\tau \int_{\mathbb{R}^d}\frac{|\vu(s,z)||\vu(s,z)|}{(|\tau-s|^{\frac{1}{\alpha}}+|y-z|)^{d+1}}dz ds\right)^{p_1}dyd\tau \right)^\frac{1}{p_1}.
\end{eqnarray*}
Extending by $0$ the function $\vu(s,z)$ when $s<0$ we obtain
\begin{eqnarray*}
\left\|\int_0^\tau \mathfrak{p}_{\tau-s}\ast\mathbb{P}\left(\mathrm{div}(\vec{u}\otimes \vec{u})\right)ds\right\|_{\mathcal{M}_{\alpha}^{p_1, \frac{d+\alpha}{\alpha-1}}}\leq C \underset{r>0}{\mathrm{sup}} \ \underset{(t,x)\in [0,+\infty[ \times \mathbb{R}^d}{\mathrm{sup}} \ \frac{1}{r^{(d+\alpha)(\frac{1}{p_1}-\frac{\alpha-1}{d+\alpha})}}\qquad \qquad\\
\times \left(\displaystyle{\iint_{\{|t-\tau|^\frac{1}{\alpha}+|x-y|<r\}}}\left( \int_{-\infty}^{+\infty} \int_{\mathbb{R}^d}\frac{|\vu(s,z)||\vu(s,z)|}{(|\tau-s|^{\frac{1}{\alpha}}+|y-z|)^{d+1}}dz ds\right)^{p_1}dyd\tau  \right)^\frac{1}{p_1},
\end{eqnarray*}
and at this point we use the definition of the parabolic Riesz potential $\mathcal{I}_{\alpha-1}$ given in (\ref{DefParaboliqueRiesz}) to write:
\begin{eqnarray*}
\left\|\int_0^\tau \mathfrak{p}_{\tau-s}\ast\mathbb{P}\left(\mathrm{div}(\vec{u}\otimes \vec{u})\right)ds\right\|_{\mathcal{M}_{\alpha}^{p_1, \frac{d+\alpha}{\alpha-1}}}&\leq & C\underset{r>0}{\mathrm{sup}} \ \underset{(t,x)\in [0,+\infty[ \times \mathbb{R}^d}{\mathrm{sup}} \ \frac{1}{r^{(d+\alpha)(\frac{1}{p_1}- \frac{\alpha-1}{d+\alpha})}}\times \\
&& \left(\displaystyle{\iint_{\{|t-\tau|^\frac{1}{\alpha}+|x-y|<r\}}}\left(\mathcal{I}_{\alpha-1}(|\vu(s,y)||\vu(s,y)|)\right)^{p_1}dyds \right)^\frac{1}{p_1}\\[3mm]
&\leq &C \|\mathcal{I}_{\alpha-1}(|\vu| |\vu|)\|_{\mathcal{M}_{\alpha}^{p_1, \frac{d+\alpha}{\alpha-1}}}.
\end{eqnarray*}
We can now use the boundedness properties in Morrey spaces of the parabolic Riesz potential $\mathcal{I}_{\alpha-1}$ given in (\ref{ContinuiteRiesz}) to obtain (since $p_1>2$):
\begin{eqnarray}
\left\|\int_0^\tau \mathfrak{p}_{\tau-s}\ast\mathbb{P}\left(\mathrm{div}(\vec{u}\otimes \vec{u})\right)ds\right\|_{\mathcal{M}_{\alpha}^{p_1, \frac{d+\alpha}{\alpha-1}}}&\leq& C\|\mathcal{I}_{\alpha-1}(|\vu| |\vu|)\|_{\mathcal{M}_{\alpha}^{p_1, \frac{d+\alpha}{\alpha-1}}}\notag\\
&\leq &C\|(|\vu| |\vu|)\|_{\mathcal{M}_{\alpha}^{\frac{p_1}{2}, \frac{d+\alpha}{2(\alpha-1)}}},\label{HolderMorreyPointFixe}
\end{eqnarray}
recall that we have $p_1> 2$ so by the H\"older inequalities in parabolic Morrey spaces (see (\ref{Holder_Morrey})) we finally have 
$$\left\|\int_0^\tau \mathfrak{p}_{\tau-s}\ast\mathbb{P}\left(\mathrm{div}(\vec{u}\otimes \vec{u})\right)ds\right\|_{\mathcal{M}_{\alpha}^{p_1, \frac{d+\alpha}{\alpha-1}}}\leq C\|\vu\|_{\mathcal{M}_{\alpha}^{p_1, \frac{d+\alpha}{\alpha-1}}}\|\vu\|_{\mathcal{M}_{\alpha}^{p_1, \frac{d+\alpha}{\alpha-1}}},$$
which is the desired estimate for the bilinear term (\ref{EstimationMorreyBilineaire}). 
\item For the external force inequality (\ref{EstimationMorreyForce}) we have:
\begin{eqnarray*}
\left\|\int_0^\tau \mathfrak{p}_{\tau-s}\ast\mathbb{P}(\vf)ds\right\|_{\mathcal{M}_{\alpha}^{p_1, \frac{d+\alpha}{\alpha-1}}}=\underset{r>0}{\mathrm{sup}} \ \underset{(t,x)\in [0,+\infty[ \times \mathbb{R}^d}{\mathrm{sup}} \ \frac{1}{r^{(d+\alpha)(\frac{1}{p_1}- \frac{\alpha-1}{d+\alpha})}}\qquad \qquad\\
\times \left(\displaystyle{\iint_{\{|t-\tau|^\frac{1}{\alpha}+|x-y|<r\}}}\left| \int_0^{\tau} \mathfrak{p}_{\tau-s}\ast\mathbb{P}(\vf)ds\right|^{p_1}dyd\tau \right)^\frac{1}{p_1},
\end{eqnarray*}
so we can write for some $0<\gamma<\alpha$:
\begin{eqnarray*}
\left\|\int_0^\tau \mathfrak{p}_{\tau-s}\ast\mathbb{P}(\vf)ds\right\|_{\mathcal{M}_{\alpha}^{p_1, \frac{d+\alpha}{\alpha-1}}}=\underset{r>0}{\mathrm{sup}} \ \underset{(t,x)\in [0,+\infty[ \times \mathbb{R}^d}{\mathrm{sup}} \ \frac{1}{r^{(d+\alpha)(\frac{1}{p_1}- \frac{\alpha-1}{d+\alpha})}}\qquad \qquad\\
\times \left(\displaystyle{\iint_{\{|t-\tau|^\frac{1}{\alpha}+|x-y|<r\}}}\left| \int_0^{\tau} \mathfrak{p}_{\tau-s}\ast\mathbb{P}((-\Delta)^{\frac{\gamma}{2}}(-\Delta)^{\frac{-\gamma}{2}}\vf)ds\right|^{p_1}dyd\tau \right)^\frac{1}{p_1}.
\end{eqnarray*}
At this point we use the following pointwise estimate for the kernel $\mathfrak{p}_{\tau-s}\ast\mathbb{P}(-\Delta)^{\frac{\gamma}{2}}(\cdot) $ 
$$|\mathfrak{p}_{t-s}\ast\mathbb{P}((-\Delta)^{\frac{\gamma}{2}}(\vpsi))(x)|\leq C\int_{\mathbb{R}^d}\frac{|\vpsi(y)|}{(|t-s|^{\frac{1}{\alpha}}+|x-y|)^{d+\gamma}}dy,$$
where $\vpsi:\mathbb{R}^d\longrightarrow \mathbb{R}^d$ is a suitable function. See a proof of this fact in the Lemma 1 of \cite{quad_non}. With this control at hand, we can now write: 
\begin{eqnarray*}
\left\|\int_0^\tau \mathfrak{p}_{\tau-s}\ast\mathbb{P}(\vf)ds\right\|_{\mathcal{M}_{\alpha}^{p, \frac{d+\alpha}{\alpha-1}}}\leq C\,\underset{r>0}{\mathrm{sup}} \ \underset{(t,x)\in [0,+\infty[ \times \mathbb{R}^d}{\mathrm{sup}} \ \frac{1}{r^{(d+\alpha)(\frac{1}{p_1}- \frac{\alpha-1}{d+\alpha})}}\qquad \qquad\\
\times \left(\displaystyle{\iint_{\{|t-\tau|^\frac{1}{\alpha}+|x-y|<r\}}}\left( \int_0^{\tau}\int_{\mathbb{R}^d}\frac{|(-\Delta)^{\frac{-\gamma}{2}}\vf(s,z)|}{(|\tau-s|^\frac{1}{\alpha}+|y-z|)^{d+\gamma}}dzds\right)^{p_1}dyd\tau \right)^\frac{1}{p_1}.
\end{eqnarray*}
We extend now the function $\vf(s,\cdot)$ by $0$ when $s<0$ and thus we can write
\begin{eqnarray*}
\left\|\int_0^\tau \mathfrak{p}_{\tau-s}\ast\mathbb{P}(\vf)ds\right\|_{\mathcal{M}_{\alpha}^{p_1, \frac{d+\alpha}{\alpha-1}}}\leq C\,\underset{r>0}{\mathrm{sup}} \ \underset{(t,x)\in [0,+\infty[ \times \mathbb{R}^d}{\mathrm{sup}} \ \frac{1}{r^{(d+\alpha)(\frac{1}{p_1}- \frac{\alpha-1}{d+\alpha})}}\qquad \qquad\\
\times \left(\displaystyle{\iint_{\{|t-\tau|^\frac{1}{\alpha}+|x-y|<r\}}}\left( \int_{-\infty}^{+\infty}\int_{\mathbb{R}^d}\frac{|(-\Delta)^{\frac{-\gamma}{2}}\vf(s,z)|}{(|\tau-s|^\frac{1}{\alpha}+|y-z|)^{d+\gamma}}dzds\right)^{p_1}dyd\tau \right)^\frac{1}{p_1}.
\end{eqnarray*}
By the definition of the parabolic Riesz potential $\mathcal{I}_{\alpha-\gamma}$ given in (\ref{DefParaboliqueRiesz}), we obtain:
\begin{eqnarray*}
\left\|\int_0^\tau \mathfrak{p}_{\tau-s}\ast\mathbb{P}(\vf)ds\right\|_{\mathcal{M}_{\alpha}^{p_1, \frac{d+\alpha}{\alpha-1}}}&\leq &C\,\underset{r>0}{\mathrm{sup}} \ \underset{(t,x)\in [0,+\infty[ \times \mathbb{R}^d}{\mathrm{sup}} \ \frac{1}{r^{(d+\alpha)(\frac{1}{p_1}- \frac{\alpha-1}{d+\alpha})}}\qquad \qquad\\
&&\times \left(\displaystyle{\iint_{\{|t-\tau|^\frac{1}{\alpha}+|x-y|<r\}}}\left(\mathcal{I}_{\alpha-\gamma}\left(|(-\Delta)^{\frac{-\gamma}{2}}\vf|\right) \right)^{p_1}dyd\tau \right)^\frac{1}{p_1}\\[3mm]
&\leq &  C \left\|\mathcal{I}_{\alpha-\gamma}\left(|(-\Delta)^{-\frac{\gamma}{2}}\vec{f}|\right)\right\|_{\mathcal{M}^{p_1,\frac{d+\alpha}{\alpha-1}}_\alpha}.
\end{eqnarray*}
Now, by the boundedness properties of the parabolic Riesz potential in parabolic Morrey spaces we have 
\begin{eqnarray}
\left\|\int_0^\tau \mathfrak{p}_{\tau-s}\ast\mathbb{P}(\vf)ds\right\|_{\mathcal{M}_{\alpha}^{p_1, \frac{d+\alpha}{\alpha-1}}}&\leq &C \left\|\mathcal{I}_{\alpha-\gamma}\left(|(-\Delta)^{-\frac{\gamma}{2}}\vec{f}|\right)\right\|_{\mathcal{M}^{p_1,\frac{d+\alpha}{\alpha-1}}_\alpha}\notag\\
&\leq &C \|(-\Delta)^{-\frac{\gamma}{2}}\vec{f}\|_{\mathcal{M}^{\mathfrak{p},\mathfrak{q}}_\alpha}=\|\vf\|_{\dot{\mathcal{W}}^{-\gamma, \mathfrak{p},\mathfrak{q}}},\label{EstimationMorreyForcePreuve}
\end{eqnarray}
where $\mathfrak{p}=\frac{(\alpha-1)p_1}{2\alpha-1-\gamma}$ and $\mathfrak{q}=\frac{d+\alpha}{2\alpha-1-\gamma}$, and this estimates gives the wished control for the external force.
\end{itemize}
We have proven the estimates (\ref{EstimationMorreyDonneeInitiale}), (\ref{EstimationMorreyBilineaire}) and (\ref{EstimationMorreyForce}), and this ends the proof of the Theorem \ref{Theorem2}.\hfill$\blacksquare$
\section{Proof of the Theorem \ref{Theorem3}}\label{Sec_Th3}
Let us start pointing out the following space inclusion: if $1<\alpha<d+2$ and $2<p_1\leq d+\alpha$, we then have
\begin{equation}\label{InclusionsMorreyMulti}
\mathcal{M}^{p_1, \frac{d+\alpha}{\alpha-1}}_\alpha([0,+\infty[\times \mathbb{R}^d)\subset \mathcal{V}_\alpha([0,+\infty[\times \mathbb{R}^d)\subset \mathcal{M}^{2, \frac{d+\alpha}{\alpha-1}}_\alpha([0,+\infty[\times \mathbb{R}^d),
\end{equation}
which is a consequence of the Fefferman-Phong inequality (see \cite{fefferman}). See a proof of these inclusions in \cite[Theorem 4]{quad_non}. These inclusions show that the functional setting of the multiplier spaces is indeed more general than the one given by the parabolic Morrey spaces and in this sense Theorem \ref{Theorem3} is a generalization of the Theorem \ref{Theorem2}.\\

Once we have specified in which sense Theorem \ref{Theorem3} is a generalization of the Theorem \ref{Theorem2}, we move on to the proof of the Theorem \ref{Theorem3} and we consider here the space $\mathcal{V}_\alpha$ as a resolution space. As in the previous results, we aim at proving the following estimates:
\begin{eqnarray}
\|\mathfrak{p}_t\ast\vu_0\|_{\mathcal{V}_\alpha}&\leq &C_1\|\vu_0\|_{\dot{B}^{-(\alpha-1), \infty}_{\infty}},\label{EstimationMultiDonneeInitiale}\\
\left\|\int_0^t \mathfrak{p}_{t-s}\ast\mathbb{P}\left(\mathrm{div}(\vec{u}\otimes \vec{u})\right)ds\right\|_{\mathcal{V}_\alpha}&\leq &C_2\|\vu\|_{\mathcal{V}_\alpha}\|\vu\|_{\mathcal{V}_\alpha},\label{EstimationMultiBilineaire}\\
\left\|\int_0^t \mathfrak{p}_{t-s}\ast\mathbb{P}(\vf)ds\right\|_{\mathcal{V}_\alpha}&\leq &C_3\|\vf\|_{\mathcal{V}_\alpha^{-1}}.\label{EstimationMultiForce}
\end{eqnarray}
\begin{itemize}
\item For the initial data, using the space inclusions (\ref{InclusionsMorreyMulti}) given above, we write:
$$\|\mathfrak{p}_t\ast\vu_0\|_{\mathcal{V}_\alpha}\leq C\|\mathfrak{p}_t\ast\vu_0\|_{\mathcal{M}^{p_1, \frac{d+\alpha}{\alpha-1}}_\alpha}\leq C\|\vu_0\|_{\dot{B}^{-(\alpha-1),\infty}_{\infty}},$$
where in the last inequality we used the control (\ref{EstimationDonneeInitiale1}). 
\item For the bilinear term we have:
$$\left|\int_0^t \mathfrak{p}_{t-s}\ast\mathbb{P}\left(\mathrm{div}(\vec{u}\otimes \vec{u})\right)ds\right|\leq C\int_{0}^t\int_{\mathbb{R}^d}\frac{|\vu(s,y)||\vu(s,y)|}{(|t-s|^\frac{1}{\alpha}+|x-y|)^{d+1}}dyds,$$
where we applied the control (\ref{fracOseen}) to estimate the kernel $\mathfrak{p}_{t-s}\ast\mathbb{P}(\mathrm{div})(\cdot)$. Setting $\vu(s,\cdot)=0$ if $s<0$, we write 
$$\left|\int_0^t \mathfrak{p}_{t-s}\ast\mathbb{P}\left(\mathrm{div}(\vec{u}\otimes \vec{u})\right)ds\right|\leq C\int_{-\infty}^{+\infty}\int_{\mathbb{R}^d}\frac{|\vu(s,y)||\vu(s,y)|}{(|t-s|^\frac{1}{\alpha}+|x-y|)^{d+1}}dyds=C\mathcal{I}_{\alpha-1}(|\vu| |\vu|)(t,x).$$
Taking the norm $\|\cdot\|_{\mathcal{V}_\alpha}$ to both sides of this estimate, we obtain
$$\left\|\int_0^t \mathfrak{p}_{t-s}\ast\mathbb{P}\left(\mathrm{div}(\vec{u}\otimes \vec{u})\right)ds\right\|_{\mathcal{V}_\alpha}\leq C\|\mathcal{I}_{\alpha-1}(|\vu| |\vu|)\|_{\mathcal{V}_\alpha}\leq C\|\vu\|_{\mathcal{V}_\alpha}\|\vu \|_{\mathcal{V}_\alpha},$$
where we applied the property (\ref{ControlRieszMultiplicateur}) and this fact proves the estimate (\ref{EstimationMultiBilineaire}).
\item For the external force we consider first the quantity 
$$\left|\int_0^t \mathfrak{p}_{t-s}\ast\mathbb{P}(\vf)ds\right|=\left|\int_0^t \mathfrak{p}_{t-s}\ast\mathbb{P}(-\Delta)^{\frac{1}{2}}((-\Delta)^{-\frac{1}{2}}\vf)ds\right|,$$
then, by the estimate (\ref{fracOseen}) (where we considered the operator $(-\Delta)^{\frac{1}{2}}$ instead of the divergence operator) we can write 
\begin{eqnarray*}
\left|\int_0^t \mathfrak{p}_{t-s}\ast\mathbb{P}(\vf)ds\right|&\leq &C\int_{0}^t\int_{\mathbb{R}^d}\frac{|(-\Delta)^{-\frac{1}{2}}\vf(s,y)|}{(|t-s|^\frac{1}{\alpha}+|x-y|)^{d+1}}dyds\\[1mm]
&\leq &C\int_{-\infty}^{+\infty}\int_{\mathbb{R}^d}\frac{|(-\Delta)^{-\frac{1}{2}}\vf(s,y)|}{(|t-s|^\frac{1}{\alpha}+|x-y|)^{d+1}}dyds,
\end{eqnarray*}
where we extended by $0$ the function $\vf(s, \cdot)$ if $s<0$. Since this is the parabolic Riesz potential $\mathcal{I}_{\alpha-1}$ given in (\ref{DefParaboliqueRiesz}), we obtain
$$\left|\int_0^t \mathfrak{p}_{t-s}\ast\mathbb{P}(\vf)ds\right|\leq C\mathcal{I}_{\alpha-1}(|(-\Delta)^{-\frac{1}{2}}\vf|)=C\mathcal{I}_{\alpha-1}\left(\sqrt{|(-\Delta)^{-\frac{1}{2}}\vf|}^2\right).$$
We take now the norm $\|\cdot\|_{\mathcal{V}_\alpha}$ to both side of the previous inequality to obtain (using the property (\ref{MonotiniciteNormeMultiplicateur}))
$$\left\|\int_0^t \mathfrak{p}_{t-s}\ast\mathbb{P}(\vf)ds\right\|_{\mathcal{V}_\alpha}\leq C\left\|\mathcal{I}_{\alpha-1}\left(\left|\sqrt{|(-\Delta)^{-\frac{1}{2}}\vf|}\right|^2\right)\right\|_{\mathcal{V}_\alpha}.$$
Recalling now that we have the equivalence  $\|\vpsi\|_{\mathcal{V}_\alpha} \simeq \sqrt{\left\| \mathcal{I}_{\alpha-1}(|\vpsi|^2) \right\|_{\mathcal{V}_\alpha}}$ given in (\ref{EquivalenceMultiplicateur}), we can write
$$\left\|\int_0^t \mathfrak{p}_{t-s}\ast\mathbb{P}(\vf)ds\right\|_{\mathcal{V}_\alpha}\leq C\left\|\sqrt{|(-\Delta)^{-\frac{1}{2}}\vf|}\right\|_{\mathcal{V}_\alpha}^2=C\|\vf\|_{\mathcal{V}_\alpha^{-1}},$$
where in the last identity we used the definition of the norm $\|\cdot\|_{\mathcal{V}_\alpha^{-1}}$ given in (\ref{EspaceMultiplicateurForce0}). We have thus proven the control (\ref{EstimationMultiForce}).
\end{itemize}
We have proven the estimates (\ref{EstimationMultiDonneeInitiale}), (\ref{EstimationMultiBilineaire}) and (\ref{EstimationMultiForce}). The proof of the Theorem \ref{Theorem3} is now ended. \hfill$\blacksquare$
\section{Counterexamples}\label{Secc_ContreExemples}
\begin{itemize}
\item In general, we do not have the inclusion $\mathcal{F}^{-\beta,p_0}_\rho([0,+\infty[\times \mathbb{R}^d)\subset \dot{\mathcal{W}}^{-\gamma, \mathfrak{p}, \mathfrak{q}} ([0,+\infty[\times \mathbb{R}^d)$ for all the admissible set of values of the parameters $\beta, p_0, \rho$ and $\gamma, \mathfrak{p}, \mathfrak{q}$ given in Theorems \ref{Theorem1} and \ref{Theorem2}.\\

In order to fix the parameters, we will consider here the dimension $d=3$ and $\alpha=\frac{3}{2}$.
\begin{itemize}
\item we consider $p_0=3$ and we have $\frac{d}{\alpha}=2<p_0<+\infty$,
\item we set $\beta=\frac{1}{3}$ so we have $0< \beta < \alpha-\frac{d}{p_0}=\frac{1}{2}$,
\item we have $\rho = 2-\tfrac{1}{\alpha}(\beta + \tfrac{d}{p_0}+1)$ and thus $\rho=\frac{4}{9}$.
\end{itemize}
For $t>0$, consider a force $\vf(t,x)=t^{-\rho}\vpsi(x)$ with $\rho=\frac{4}{9}$ and where $\vpsi:\mathbb{R}^3\longrightarrow \mathbb{R}^3$ is a function such that $(-\Delta)^{-\frac{\beta}{2}}\vpsi\in L^{3}(\mathbb{R}^3)$.\\

With these conditions, we have $\vf\in \mathcal{F}^{-\beta,p_0}_\rho([0,+\infty[\times \mathbb{R}^3)$, indeed 
$$\|\vec{f}\|_{\mathcal{F}^{-\beta,p_0}_\rho}=\underset{t>0}{\mathrm{sup}} \ \tau^\rho \; \|t^{-\rho}(-\Delta)^{\frac{-\beta}{2}}\vpsi(\cdot)\|_{L^{p_0}}\leq \|(-\Delta)^{\frac{-\beta}{2}}\vpsi(\cdot)\|_{L^{p_0}}<+\infty.$$
Thus, from a initial data $\vu_0\in \dot{B}^{-\frac{1}{2},\infty}_{\infty}(\mathbb{R}^3)$ and for this particular external force $\vf$, by the Theorem \ref{Theorem1}, we can obtain a mild solution for the fractional Navier-Stokes equations in the space $L^{\infty}_{\frac{3}{2}}(\mathbb{R}^3)$.\\

However, in general we do not have $\vf\in \dot{\mathcal{W}}^{-\gamma, \mathfrak{p}, \mathfrak{q}}([0,+\infty[\times \mathbb{R}^d)$. Indeed, if we set 
\begin{itemize}
\item we set $p_1=2.9$, we have $2<p_1<\frac{\alpha}{\alpha-1}=3$,
\item if we have $\gamma=1.4$, then we have $2\alpha-1-(\alpha-1)p_1<\gamma<\alpha$,
\item we set $\mathfrak{p}=\frac{(\alpha-1)p_1}{2\alpha-1-\gamma}\sim 2.41666>\frac{9}{4}=2.25$.
\end{itemize}
Thus all the conditions given in Theorem \ref{Theorem2} for the indexes that define the space $\dot{\mathcal{W}}^{-\gamma, \mathfrak{p}, \mathfrak{q}}$ are fulfilled. But if we consider only the time variable (\emph{i.e.} if we assume that we have $(-\Delta)^{\frac{-\gamma}{2}}\vpsi\in L^{\mathfrak{p}}(\mathbb{R}^3)$) we need to study the integrability of the quantity $s^{-\rho \mathfrak{p}}$, and we have $\rho \mathfrak{p}>1$ (since $\rho=\frac{4}{9}$ and $\mathfrak{p}>\frac{9}{4}$) and thus this function is not integrable on intervals that contain the origin (it is not locally integrable), so we obtain 
$$\|\vf\|_{\dot{\mathcal{W}}^{-\gamma, \mathfrak{p}, \mathfrak{q}}}=+\infty.$$
From this we see that we do not have in general and for \emph{all} the admissible parameters $\beta, p_0, \rho$ and $\gamma, \mathfrak{p}, \mathfrak{q}$ the inclusion $\mathcal{F}^{-\beta,p_0}_\rho([0,+\infty[\times \mathbb{R}^d)\subset \dot{\mathcal{W}}^{-\gamma, \mathfrak{p}, \mathfrak{q}} ([0,+\infty[\times \mathbb{R}^d)$.

An interesting consequence of this fact is that, for an initial data $\vu_0\in \dot{B}^{-\frac{1}{2},\infty}_{\infty}(\mathbb{R}^3)$ and for this external force $\vf\in \mathcal{F}^{-\beta,p_0}_\rho$ (\emph{i.e.} for this set of parameters $\beta, p_0, \rho$), we can not close the fixed-point argument in the parabolic Morrey space $\mathcal{M}^{p_1, \frac{d+\alpha}{\alpha-1}}_\alpha([0,+\infty[\times \mathbb{R}^3)$ with the space $\dot{\mathcal{W}}^{-\gamma, \mathfrak{p}, \mathfrak{q}} ([0,+\infty[\times \mathbb{R}^d)$ for the forces. \\

For the sake of completeness, we give now some conditions on the parameters $\beta, p_0, \rho$ and $\gamma, \mathfrak{p}, \mathfrak{q}$ (without any claim of optimality) that ensure the space embedding $\mathcal{F}^{-\beta,p_0}_\rho([0,+\infty[\times \mathbb{R}^d)\subset \dot{\mathcal{W}}^{-\gamma, \mathfrak{p}, \mathfrak{q}} ([0,+\infty[\times \mathbb{R}^d)$. In this sense we have the following result:
\begin{Lemme}
Consider the spaces $\mathcal{F}^{-\beta,p_0}_\rho([0,+\infty[\times \mathbb{R}^d)$ and $\dot{\mathcal{W}}^{-\gamma, \mathfrak{p}, \mathfrak{q}} ([0,+\infty[\times \mathbb{R}^d)$  as defined in (\ref{Def_NormeF1}) and (\ref{Definition_MorreySobolev}), respectively. If we have $\rho \mathfrak{p}<1$, $\mathfrak{p}\leq p_0$ and $\beta=\gamma$, then we have the space inclusion
$$\mathcal{F}^{-\beta,p_0}_\rho([0,+\infty[\times \mathbb{R}^d)\subset \dot{\mathcal{W}}^{-\gamma, \mathfrak{p}, \mathfrak{q}}([0,+\infty[\times \mathbb{R}^d).$$
\end{Lemme}
\noindent{\bf Proof.} We start by considering, for some fixed $t>0$ and $x\in \mathbb{R}^3$, the quantity 
$$I=\frac{1}{r^{(d+\alpha)(1-\frac{\mathfrak{p}}{\mathfrak{q}})}} \displaystyle{\int_{\{|t-s|<r^\alpha\}}}\displaystyle{\int_{\{|x-y|<r\}}}|(-\Delta)^{-\frac{\beta}{2}}\vec{f}(s,y)|^\mathfrak{p}dyds.$$
Applying the Hölder inequality with $1=\frac{\mathfrak{p}}{p_0}+\frac{p_0-\mathfrak{p}}{p_0}$ on the integral with respect to the space variable (recall that by hypothesis we have $1\leq \frac{p_0}{\mathfrak{p}}$) we have the estimate
\begin{eqnarray*}
I&\leq &\frac{1}{r^{(d+\alpha)(1-\frac{\mathfrak{p}} {\mathfrak{q}})}} \displaystyle{\int_{|t-s|<r^\alpha}} \left(\displaystyle{\int_{|x-y|<r}}|(-\Delta)^{-\frac{\beta}{2}}\vec{f}(s,x)|^{p_0}dy\right)^\frac{\mathfrak{p}}{p_0} |B_r|^{\frac{p_0-\mathfrak{p}}{p_0}}ds\\
&\leq &\frac{Cr^{d(\frac{p_0-\mathfrak{p}}{p_0})}}{r^{(d+\alpha)(1-\frac{\mathfrak{p}} {\mathfrak{q}})}} \displaystyle{\int_{|t-s|<r^\alpha}} \left(\displaystyle{\int_{|x-y|<r}}|(-\Delta)^{-\frac{\beta}{2}}\vec{f}(s,x)|^{p_0}dy\right)^\frac{\mathfrak{p}}{p_0}ds,
\end{eqnarray*}
and taking the integral in space on the entire domain $\mathbb{R}^d$, the $L^{p_0}(\mathbb{R}^d)$ norm appears, allowing us to obtain
\begin{eqnarray*}
I&\leq &\frac{Cr^{d(\frac{p_0-\mathfrak{p}}{p_0})}}{r^{(d+\alpha)(1-\frac{\mathfrak{p}} {\mathfrak{q}})}} \displaystyle{\int_{|t-s|<r^\alpha}} \left(\|(-\Delta)^{-\frac{\beta}{2}}\vec{f}(s,\cdot)\|_{L^{p_0}}\right)^{\mathfrak{p}}ds\\
&\leq &\frac{Cr^{d(\frac{p_0-\mathfrak{p}}{p_0})}}{r^{(d+\alpha)(1-\frac{\mathfrak{p}} {\mathfrak{q}})}} \displaystyle{\int_{|t-s|<r^\alpha}} \left(|s|^\rho\|(-\Delta)^{-\frac{\beta}{2}}\vec{f}(s,\cdot)\|_{L^{p_0}}\right)^{\mathfrak{p}}|s|^{-\rho\mathfrak{p}}ds.
\end{eqnarray*}
Recalling that $\|\vec{f}\|_{\mathcal{F}^{-\beta,p_0}_\rho}  = \underset{s>0}{\mathrm{sup}} \ s^\rho\|(-\Delta)^{\frac{-\beta}{2}}\vf(s, \cdot)\|_{L^{p_0}}$ we can write 
$$I\leq \|\vf\|_{\mathcal{F}^{-\beta,p_0}_\rho}^{\mathfrak{p}}\frac{Cr^{d(\frac{p_0-\mathfrak{p}}{p_0})}}{r^{(d+\alpha)(1-\frac{\mathfrak{p}} {\mathfrak{q}})}} \displaystyle{\int_{|t-s|<r^\alpha}} |s|^{-\rho\mathfrak{p}}ds.$$
Now, since by hypothesis we have $\rho\mathfrak{p}<1$, the previous integral is finite and we have
\begin{equation}\label{PuissanceRInclusion}
I\leq \|\vf\|_{\mathcal{F}^{-\beta,p_0}_\rho}^{\mathfrak{p}}\frac{Cr^{d(\frac{p_0-\mathfrak{p}}{p_0})}}{r^{(d+\alpha)(1-\frac{\mathfrak{p}} {\mathfrak{q}})}}r^{\alpha(1-\rho\mathfrak{p})}= C\|\vf\|_{\mathcal{F}^{-\beta,p_0}_\rho}^{\mathfrak{p}} r^{d(\frac{p_0-\mathfrak{p}}{p_0})+\alpha(1-\rho\mathfrak{p})-(d+\alpha)(1-\frac{\mathfrak{p}} {\mathfrak{q}})}.
\end{equation}
We study now the power of the parameter $r$ and we have 
\begin{eqnarray*}
d(\tfrac{p_0-\mathfrak{p}}{p_0})+\alpha(1-\rho\mathfrak{p})-(d+\alpha)(1-\tfrac{\mathfrak{p}} {\mathfrak{q}})&=&d-\tfrac{d\mathfrak{p}}{p_0}+\alpha-\alpha\rho\mathfrak{p} - d+d\tfrac{\mathfrak{p}} {\mathfrak{q}}-\alpha+\alpha\tfrac{\mathfrak{p}} {\mathfrak{q}}\\
&=&\mathfrak{p}(\tfrac{d+\alpha}{\mathfrak{q}}-\tfrac{d}{p_0}-\alpha\rho),
\end{eqnarray*}
but since we have $\alpha\rho= 2\alpha-(\beta + \tfrac{d}{p_0}+1)$ by (\ref{Def_NormeF1}) and since we have $\mathfrak{q}=\frac{d+\alpha}{2\alpha-1-\gamma}$, we can write
\begin{eqnarray*}
d(\tfrac{p_0-\mathfrak{p}}{p_0})+\alpha(1-\rho\mathfrak{p})-(d+\alpha)(1-\tfrac{\mathfrak{p}} {\mathfrak{q}})&=&\mathfrak{p}\left((2\alpha-1-\gamma)-\tfrac{d}{p_0}-2\alpha+(\beta + \tfrac{d}{p_0}+1)\right)\\
&=&\mathfrak{p}(\beta-\gamma),
\end{eqnarray*}
but since $\beta=\gamma$ by hypothesis, we finally have that the power of the parameter $r$ in (\ref{PuissanceRInclusion}) is null, so we obtain the uniform (in $t>0$, $x\in \mathbb{R}^3$ and $r>0$) estimate $I\leq  C\|\vf\|_{\mathcal{F}^{-\beta,p_0}_\rho}^{\mathfrak{p}}$, from which we deduce the inequality 
$$\|\vf\|_{\dot{\mathcal{W}}^{-\gamma, \mathfrak{p}, \mathfrak{q}}}\leq C\|\vf\|_{\mathcal{F}^{-\beta,p_0}_\rho},$$
and the lemma is proven. \hfill $\blacksquare$\\
\item We study now the inclusion $\dot{\mathcal{W}}^{-\gamma, \mathfrak{p}, \mathfrak{q}} ([0,+\infty[\times \mathbb{R}^d)\subset \mathcal{V}^{-1}_\alpha([0,+\infty[\times \mathbb{R}^d)$. As above, we will consider here the dimension $d=3$, $\alpha=\frac{3}{2}$ and we set:
\begin{itemize}
\item $p_1=\frac{5}{2}$, so we have  $2<p_1<\frac{\alpha}{\alpha-1}=3$,
\item $\gamma=\frac{7}{5}$, and we get $2\alpha-1-(\alpha-1)p_1=\frac{3}{4}<\gamma<\alpha$,
\item then we set $\mathfrak{p}=\frac{(\alpha-1)p_1}{2\alpha-1-\gamma}=\frac{25}{12}$,
\item and $\mathfrak{q}=\frac{d+\alpha}{2\alpha-1-\gamma}=\frac{15}{2}$.
\end{itemize}
With these parameters, we consider now a force $\vec{g}:[0,+\infty[\times \mathbb{R}^3\longrightarrow \mathbb{R}^3$ given by the following expression
\begin{equation}\label{DefForceCtreExTh23}
\vec{g}(t,x)=|t|^{-\frac{2}{5}}\cos(e_0 \cdot x)v_0,
\end{equation}
where $e_0$ and $v_0$ are vectors in $\mathbb{R}^3$ such that $|e_0|=|v_0|=1$. Note that for $s\in \mathbb{R}$ and since $e_0\neq0$, we can define the fractional laplacian of $\vg$ as follows
$$((-\Delta)^\frac{s}{2}\vec{g}(t,\cdot))^{\widehat{\quad}}(\xi) = C|t|^{-\frac{2}{5}}|\xi|^s\frac{\delta_{-e_0}(\xi)+\delta_{e_0}(\xi)}{2}\vec{v_0}.$$
Using the fact that $|e_0|=1$ and the properties of Dirac masses, we obtain
$$((-\Delta)^\frac{s}{2}\vec{g}(t,\cdot))^{\widehat{\quad}}(\xi) = C|t|^{-\frac{2}{5}}\frac{|-e_0|^s\delta_{-e_0}(\xi)+|e_0|^s\delta_{e_0}(\xi)}{2}\vec{v_0} = C|t|^{-\frac{2}{5}}\frac{\delta_{-e_0}(\xi)+\delta_{e_0}(\xi)}{2}\vec{v_0},$$
and we obtain $((-\Delta)^\frac{s}{2}\vec{g}(t,\cdot))^{\widehat{}}(\xi) =\vec{g}(t,.)^{\widehat{}}(\xi)$ in other words, we have for all $s\in \mathbb{R}$
$$(-\Delta)^\frac{s}{2}\vec{g}(t,x)=\vec{g}(t,x).$$
With this remark at hand we thus have $\vg\in\dot{\mathcal{W}}^{-\frac{7}{5}, \frac{25}{12}, \frac{15}{2}}([0,+\infty[\times \mathbb{R}^3)$, indeed we write 
$$|(-\Delta)^{-\frac{7}{5}}\vg(t,x)|=| |t|^{-\frac{2}{5}}(-\Delta)^{-\frac{7}{5}}\cos(e_0 \cdot x)v_0|,$$
and since $|e_0|=|v_0|=1$, we obtain 
\begin{equation}\label{EstimationContrexTH23}
|(-\Delta)^{-\frac{7}{5}}\vg(t,x)|=|t|^{-\frac{2}{5}}|\cos(e_0 \cdot x)|\leq |t|^{-\frac{2}{5}}.
\end{equation}
Now, with this estimate we write
\begin{eqnarray*}
I=\frac{1}{r^{(d+\alpha)(\frac{1}{\mathfrak{p}}-\frac{1}{\mathfrak{q}})}} \left(\int_{\{|t-s|<r^\alpha\}}\int_{\{|x-y|<r\}}|(-\Delta)^{-\frac{7}{5}}\vg(s,y)|^\mathfrak{p}dyds\right)^{\frac{1}{\mathfrak{p}}}\\
=\frac{1}{r^{\frac{9}{2}(\frac{12}{25}-\frac{2}{15})}}\left(\int_{\{|t-s|<r^{\frac{3}{2}}\}}\int_{\{|x-y|<r\}} |(-\Delta)^{-\frac{7}{5}}\vg(s,y)|^\frac{25}{12}dyds\right)^{\frac{12}{25}},
\end{eqnarray*}
and we obtain
\begin{eqnarray*}
I&\leq &C\frac{1}{r^{\frac{9}{2}(\frac{12}{25}-\frac{2}{15})}}\left(\int_{\{|t-s|<r^{\frac{3}{2}}\}}\int_{\{|x-y|<r\}} \left(|t|^{-\frac{2}{5}}\right)^\frac{25}{12}dyds\right)^{\frac{12}{25}}\\
&\leq &C\frac{1}{r^{\frac{39}{25}}}\left(\int_{\{|t-s|<r^{\frac{3}{2}}\}}|s|^{-\frac{5}{6}}ds \int_{\{|x-y|<r\}}dy\right)^{\frac{12}{25}}.
\end{eqnarray*}
By evaluating the integrals above we have uniform control 
$$I\leq C\frac{1}{r^{\frac{39}{25}}}\times r^{\frac{39}{25}}\leq C,$$
from which we easily deduce that $\vg\in \dot{\mathcal{W}}^{-\frac{7}{5}, \frac{25}{12}, \frac{15}{2}}([0,+\infty[\times \mathbb{R}^3)$.\\

Now, by the definition of the space $\mathcal{V}^{-1}_\alpha$ given in (\ref{EspaceMultiplicateurForce0}) if we have $\vg\in \mathcal{V}^{-1}_\alpha$ we should have $\sqrt{|(-\Delta)^{-\frac{1}{2}}\vg|}\in \mathcal{V}_\alpha$, thus by the inclusions (\ref{InclusionsMorreyMulti}) we have $\sqrt{|(-\Delta)^{-\frac{1}{2}}\vg|}\in\mathcal{M}^{2,\frac{d+\alpha}{\alpha-1}}=\mathcal{M}^{2,9}$, from which we deduce that we should have $|(-\Delta)^{-\frac{1}{2}}\vg|\in\mathcal{M}^{1,\frac{9}{2}}$. However, by the same computations that leaded us to the expression (\ref{EstimationContrexTH23}), we have 
$$|(-\Delta)^{-\frac{1}{2}}\vg(t,x)|=|t|^{-\frac{2}{5}}|\cos(e_0 \cdot x)|,$$
but this function does not belong to the Morrey space $\mathcal{M}^{1,\frac{9}{2}}$. Indeed, if $|(-\Delta)^{-\frac{1}{2}}\vg(t,x)|\in \mathcal{M}^{1,\frac{9}{2}}$ we should have
$$\|(-\Delta)^{-\frac{1}{2}}\vg\|_{\mathcal{M}^{1,\frac{9}{2}}}=\underset{r>0}{\mathrm{sup}} \ \underset{(t,x)\in [0,+\infty[ \times \mathbb{R}^d}{\mathrm{sup}} \ \frac{1}{r^{\frac{7}{2}}} \displaystyle{\int_{\{|t-s|<r^\frac{3}{2}\}}}\displaystyle{\int_{\{|x-y|<r\}}}|(-\Delta)^{-\frac{1}{2}}\vec{g}(s,y)|dyds<+\infty,$$
but 
$$\|(-\Delta)^{-\frac{1}{2}}\vg\|_{\mathcal{M}^{1,\frac{9}{2}}}\geq \underset{r>0}{\mathrm{sup}} \ \frac{1}{r^{\frac{7}{2}}} \int_{\{|s|<r^\frac{3}{2}\}}\int_{\{|y|<r\}}  |s|^{-\frac{2}{5}}|\cos(e_0 \cdot y)|dyds,$$
by looking of the contribution of the function $|\cos(e_0 \cdot y)|$ over balls $B(0,r)$ of radius $r>0$, we can obtain (for some constant $C\gg 1$)
$$\|(-\Delta)^{-\frac{1}{2}}\vg\|_{\mathcal{M}^{1,\frac{9}{2}}}\geq \frac{1}{C}\ \underset{r>0}{\mathrm{sup}} \ \frac{1}{r^{\frac{7}{2}}} \int_{\{|s|<r^\frac{3}{2}\}}|s|^{-\frac{2}{5}}ds \ |B(0,r)|,$$
thus, we see that 
$$\int_{\{|s|<r^\frac{3}{2}\}}|s|^{-\frac{2}{5}}ds \ |B(0,r)|=C'r^{\frac{3}{2}(1-\frac{2}{5})}r^3=C'r^{\frac{39}{10}},$$
and the quantity $r^{\frac{39}{10}}$ can not be compensated by the weight $r^{\frac{7}{2}}$ for all $r>0$ and this shows that $|(-\Delta)^{-\frac{1}{2}}\vg(t,x)|\notin \mathcal{M}^{1,\frac{9}{2}}$.\\

From this, we can see that there is not a simple relationship between the spaces
$$\dot{\mathcal{W}}^{-\gamma, \mathfrak{p}, \mathfrak{q}}([0,+\infty[\times \mathbb{R}^d)\qquad \mbox{and}\qquad \mathcal{V}^{-1}_\alpha([0,+\infty[\times \mathbb{R}^d),$$
for \emph{all} the admissible parameters $\gamma, \mathfrak{p}, \mathfrak{q}$ given in Theorem \ref{Theorem2}.\\

Thus, just as in the previous point, if $\vu_0\in \dot{B}^{-\frac{1}{2},\infty}_{\infty}(\mathbb{R}^3)$ and for this particular force (\emph{i.e.} for the parameters $\gamma, \mathfrak{p}$ and $\mathfrak{q}$ given above) $\vg\in\dot{\mathcal{W}}^{-\gamma, \mathfrak{p}, \mathfrak{q}}([0,+\infty[\times \mathbb{R}^d)$, we can construct by the Theorem \ref{Theorem2} a mild solution for the fractional Navier-Stokes equation in the parabolic Morrey space $\mathcal{M}^{p_1, \frac{d+\alpha}{\alpha-1}}_\alpha([0,+\infty[\times \mathbb{R}^3)$ but it is not possible to close the fixed-point argument for the data $(\vu_0, \vg)$ in the space $\mathcal{V}_\alpha([0,+\infty[\times \mathbb{R}^3)$ used in Theorem \ref{Theorem3} with the space $\mathcal{V}^{-1}_\alpha([0,+\infty[\times \mathbb{R}^3)$ for the force $\vg$.\\

As before and for the sake of completeness, we give here one simple criterion for the space inclusion
$\dot{\mathcal{W}}^{-\gamma, \mathfrak{p}, \mathfrak{q}} \subset \mathcal{V}^{-1}_\alpha$.
\begin{Lemme}
For $d\geq 2$ and for $1<\alpha<2$, if $\gamma=1$, $\mathfrak{p}=\frac{(\alpha-1)p_1}{2\alpha-1-\gamma}=\frac{p_1}{2}$ and $\mathfrak{q}=\frac{d+\alpha}{2\alpha-1-\gamma}=\frac{d+\alpha}{2\alpha-2}$, then we have the embedding 
$$\dot{\mathcal{W}}^{-\gamma, \mathfrak{p}, \mathfrak{q}}([0,+\infty[\times \mathbb{R}^3) \subset \mathcal{V}^{-1}_\alpha([0,+\infty[\times \mathbb{R}^3).$$
\end{Lemme}
{\bf Proof.} Recall that by (\ref{EspaceMultiplicateurForce0}) we have 
$$\|\vec{f}\|_{\mathcal{V}_\alpha^{-1}}=\left\|\sqrt{|(-\Delta)^{-\frac{1}{2}}\vec{f}|}\right\|_{\mathcal{V}_\alpha}^2,$$
recall also that by (\ref{EquivalenceMultiplicateur}) we have the equivalence 
$$\|\vec{f}\|_{\mathcal{V}_\alpha^{-1}}=\left\|\sqrt{|(-\Delta)^{-\frac{1}{2}}\vec{f}|}\right\|_{\mathcal{V}_\alpha}^2\simeq \left\|\mathcal{I}_{\alpha-1}\left(|(-\Delta)^{-\frac{1}{2}}\vec{f}|\right)\right\|_{\mathcal{V}_\alpha},$$
now by the Fefferman-Phong inequalities (\ref{InclusionsMorreyMulti}) we can write 
$$\left\|\mathcal{I}_{\alpha-1}\left(|(-\Delta)^{-\frac{1}{2}}\vec{f}|\right)\right\|_{\mathcal{V}_\alpha}\leq C\left\|\mathcal{I}_{\alpha-1}\left(|(-\Delta)^{-\frac{1}{2}}\vec{f}|\right)\right\|_{\mathcal{M}^{p_1, \frac{d+\alpha}{\alpha-1}}_\alpha}.$$
At this point we use the boundedness of the Riesz potentials in Morrey spaces given in (\ref{ContinuiteRiesz}) to obtain
$$\left\|\mathcal{I}_{\alpha-1}\left(|(-\Delta)^{-\frac{1}{2}}\vec{f}|\right)\right\|_{\mathcal{M}^{p_1, \frac{d+\alpha}{\alpha-1}}_\alpha}\leq C\left\||(-\Delta)^{-\frac{1}{2}}\vec{f}|\right\|_{\mathcal{M}^{\mathfrak{p},\mathfrak{q}}_\alpha},$$
with $\mathfrak{p}=\frac{(\alpha-1)p_1}{2(\alpha-1)}=\frac{p_1}{2}$ and $\mathfrak{q}=\frac{d+\alpha}{2\alpha-2}$. We can thus write by (\ref{Definition_MorreySobolev})
$$\|\vec{f}\|_{\mathcal{V}_\alpha^{-1}}\leq C\|\vf\|_{\dot{\mathcal{W}}^{-1,\mathfrak{p}, \mathfrak{q}}},$$
from which we deduce the wished space inclusion. \hfill $\blacksquare$
\end{itemize}



\end{document}